\documentclass[12pt,a4paper,final]{article}
\usepackage[margin=2.54cm]{geometry}%

\usepackage{epsfig,amsmath,amsfonts,bef_alex,amssymb}
\usepackage{eucal,bm,color,tikz,scrtime,todonotes}
\numberwithin{equation}{section}
\numberwithin{figure}{section}
\usepackage[utf8]{inputenc}

\usepackage[notref,notcite]{showkeys}

\newcommand{\mfT}{\mathfrak T}

\newcommand{\sign}{\mathop{\mafo{sign}}}
\newcommand{\Sign}{\mathop{\mafo{Sign}}}
\newcommand{\DISS}{\mafo{Var}_{\|\cdot\|}}

\newcommand{\COMM}[1]{}

\begin{document}

\author{Alexander Mielke\thanks{Weierstra\ss{}-Institut f\"ur Angewandte Analysis und
    Stochastik, Anton-Wilhelm-Amo-Stra\ss{}e 39, 10117 Berlin and Humboldt
    Universit\"at zu Berlin, Germany.}}

\title{Relating a rate-independent system and \\ a gradient system for
  the case of\\ one-homogeneous potentials\thanks{Research
    partially supported by DFG via SFB\,910 (project no.\ 163436311), subproject A5.}}

\date{30. September 2020}
 
\maketitle

\begin{center} 
\itshape\large In memory of Genevi\`eve Raugel,\\
a brilliant mathematician and a caring friend
\end{center}


\section{Introduction}
\label{se:Intro}

In this paper we consider two types of generalized gradient systems
$(X,\calF,\calR)$ in the
sense of \cite{Miel16EGCG}, which are both given in terms of a Hilbert space $X$,
an energy functional $\calF$, and a dissipation structure $\calR$, such that the
induced evolution equation takes the form
$0 \in \pl\calR(\dot u(t)) + \pl\calF(t,u(t))$.  The two systems are
\begin{description}
\item[the Hilbert-space gradient system (GS)] $(X,\calJ,\frac12\|\cdot\|^2)$ and

\item[the energetic rate-independent system (ERIS)] $(X,\calE,\|\cdot\|)$ with
  $\calE(t,u)= t \calJ(u)$. 
\end{description}
Here $X$ is a Hilbert space with duality pairing $\langle \cdot,\cdot\rangle$,
norm $\|\cdot\|$, and Riesz isomorphism $E:X\to X^*$.  The link between these
two system arises from the fact that we assume that $\calJ:X\to [0,\infty]$ is
positively homogeneous of degree 1, i.e.\ $\calJ(\lambda u)=\lambda \calJ(u)$
for all $\lambda> 0$ and $u\in X$. Moreover, $\calJ$ is convex, lower
semicontinuous, and has a dense domain
dom$(\calJ):=\bigset{u\in X}{\calJ(u)<\infty} $.

The gradient-flow equation for the gradient system
$(X,\calJ,\frac12\|\cdot\|^2)$ takes the form
\begin{equation}
  \label{eq:GFE}
  0 \in E w'(s) + \pl\calJ(w(s)) \quad \text{for a.a. }s>0, \qquad w(0)=u^0.  
\end{equation}
We continue to use the letter $s\geq 0$ for the time in the gradient-flow
equation, while $t\geq 0$ will be reserved for the time in the ERIS.
On the formal level, the evolution equation induced by the ERIS
$(X,\calE,\|\cdot\|)$ can be written in the analogous form 
\begin{equation}
  \label{eq:DNE}
  0 \in \mathop{\mafo{Sign}}(\dot u(t)) + \pl_u \calE(t,u(t)) \quad \text{for
    a.a.\ } t>0, \qquad u(0)=u^0 ,
\end{equation}
where $\mathop{\mafo{Sign}}(v)\subset X^*$ denotes the Hilbert-space signum
function, which is the convex subdifferential
of the norm $\|\cdot\|$, namely
\[
\mathop{\mafo{Sign}}(v) = \pl(\|\cdot\|)(v) = \left\{ \ba{cl} 
\bigset{\xi\in X^* }{\|\xi\|_*\leq 1} & \text{for } v=0, \\[0.3em]
\ds \frac1{\|v\|} Ev & \text{for } v\neq 0. \ea\right.
\]

Recalling $\calE(t,u)=t\calJ(u)$, we find $\pl_u \calE(t,u)= t \pl\calJ(u)$
and see on the formal level that the gradient-flow equation \eqref{eq:GFE}
and the rate-independent evolution \eqref{eq:DNE} are equivalent up to a
time reparametrization. Indeed, assuming that  $u:{[0,\infty[} \to X$ is a sufficiently
smooth solution of rate-independent evolution \eqref{eq:DNE}, we define a
reparametrization $s=S(t)$ via 
\[
S(t)= \int_0^t \tau \| \dot u(\tau)\| \dd \tau
\]
and assume further that we can invert the relation to obtain $t=T(s)$. 
Then, the chain rule shows that $w$ defined via $w(s)=u(T(s))$ is a solution of
the gradient-flow equation \eqref{eq:GFE}.
Vice versa,  if a sufficiently smooth solution $w$ of \eqref{eq:GFE} 
is given, we define the
reparametrization $t=T(s)$ via 
\begin{equation}
  \label{eq:t=T(s)}
  T(s)= \frac1{\|w'(s)\|} 
\end{equation}
and assume that the inversion $s=S(t)$ exists, then $u(t)=w(S(t))$ solves
\eqref{eq:DNE}. 

Because of the 1-homogeneity of $\calJ$ and the simple structure of the
dissipation in terms of the norm $\|\cdot\|$, both systems have a scaling
invariance. For all $\lambda>0$ we have the implications:
\begin{subequations}
  \label{eq:Scaling}
\begin{align}
w \text{ solves \eqref{eq:GFE}} \quad &\Longrightarrow \quad  w_\lambda: s
\mapsto \frac1\lambda\,w(\lambda s) \text{ solves \eqref{eq:GFE}},
\\
u \text{ solves \eqref{eq:DNE}} \quad &\Longrightarrow \quad  \wt u_\lambda: s
\mapsto\ \ u (\lambda s) \ \  \text{ solves \eqref{eq:DNE}}.
\end{align}
\end{subequations}
As for linear equations $w'=-Lw$, where the existence of an eigenpair
$(\varphi,\lambda)$ of $L$, i.e.\ $L\varphi =\lambda \varphi$ leads to the
explicit solutions $w(s)= c\ee^{-\lambda s} \varphi$, we obtain explicit
solutions from nontrivial solutions of the relation $\pl^0\calJ(\psi)=\lambda
E\psi$, where $\pl^0\calJ(u)$ denotes the unique element in $\pl\calJ(u)$ with
minimal norm $\|u\|$. For all $\rho>0$ we find that
\begin{subequations}
 \label{eq:PiecewiseSoln}
 \begin{align}
  \label{eq:PiecewiseSoln.w}
 &w(s)=\max\{\rho{-}\lambda s, 0\}\, \psi \ \text{ solves \eqref{eq:GFE}, \qquad
   and} \\ 
   \label{eq:PiecewiseSoln.u}
 &u(t)= \bm1_{[0,t_*]}(t) \,\rho\psi \ \text{ solves \eqref{eq:DNE} \quad (in the
  sense of \eqref{eq:def.ES}),} 
 \end{align}
\end{subequations}
where $t_*=1/\|\lambda \psi\|$ and $\bm1_A(t)=1$ for $t\in A$ and $0$
otherwise.\medskip  

The purpose of this work is to make these formal observations rigorous, thus
relating the two generalized gradient systems in an mathematically precise
way. 

To indicate one of the difficulties, we observe that for \eqref{eq:t=T(s)} the
mapping $s\mapsto \|w'(s)\|$ should be monotonous, which is a standard feature
for Hilbert-space gradient flows with convex potentials (cf.\
\cite[Thm.\,3.1(6)]{Brez73OMMS}), but as teh solution in
\eqref{eq:PiecewiseSoln.w} shows, we cannot expect strict monotonicity.
This is indeed related to the fact that the solutions of the rate-independent
evolution \eqref{eq:DNE} are not even continuous, see
\eqref{eq:PiecewiseSoln.u} for an example. Hence, it is necessary to
replace the differential formulation \eqref{eq:DNE} by a derivative-free one,
which is available for ERIS. We refer to \cite{Miel08CIME, MieRou15RIST} for
different solution concepts of rate-independent systems. 

For our purposes, the
concept of \emph{energetic solutions} for the ERIS as introduced in
\cite{MiThLe02VFRI, Miel05ERIS} will be appropriate as it allows for jumps.  We
call $u:{[0,\infty[} \to X$ an energetic solution for $(X,\calE,\|\cdot\|)$, if
$t \mapsto \pl_t \calE(t,u(t))$ lies in $\rmL^1_\mafo{loc}({[0,\infty[})$ and
for all $t\geq 0$ we have the global stability (S) and the energy balance (E):
\begin{equation}
  \label{eq:def.ES}
  \begin{aligned}
\text{(S)}\qquad &\forall\,\wt u\in X:\quad \calE(t,\wt u)\leq
\calE(t,u(t))+ \|\wt u{-}u(t)\|; \\
\text{(E)}\qquad & \calE(t,u(t)) + \text{Var}_{\|\cdot\|}(u,[0,t]) =
\calE(0,u(0)) + \int_0^t \pl_s \calE(s,u(s))\dd s,  
\end{aligned}
\end{equation}
where $\text{Var}_{\|\cdot\|}(u,[s,t])= \sup\set{\sum_{j=1}^N \|
  u(t_j){-}u(t_{j-1}) \|}{ N\in \N,\ s\leq t_0<t_1< \cdots < t_N\leq
  t}$. 

Because of the convexity of $u\mapsto \calE(t,u)=t\calJ(u)$ it is obvious to
construct approximate solutions via the minimizing-movement scheme, i.e.\ by
choosing a time step $h>0$ and defining $u^{k}_h$ as the minimizer of
$u \mapsto \| u{-}u^{k-1}_h\| + \calE(k h, u)$. Under the additional assumption
that $\calJ$ has compact sublevels in $X$, it is then standard (cf.\
\cite{Miel05ERIS,MieRou15RIST}) to show that solutions can be obtained as
accumulation points of these approximations, see also Proposition
\ref{pr:Exist.Cmpt}. However, without these compactness assumption the
existence of solutions is largely open, except for the case that
$\calE(t,\cdot)$ is quadratic.  Even worse, uniqueness can only be shown in
situations where $\calE(t,\cdot)$ has a Lipschitz-continuous second
derivative, see \cite{MieThe04RIHM, BrKrSc04UEQI, MieRos07EURC}. Thus, it is
surprising that the ERIS with $\calE(t,u)=t\calJ(u)$ as considered here provides
a model class, where we can show both, (i) existence of solutions
without assuming compactness and (ii) uniqueness of solutions.

Here uniqueness holds up to the choice of the jump behavior. Because of
the finiteness of the variation $\text{Var}_{\|\cdot\|}(u,[0,t])$ it is clear
that at all times $t\geq 0$ the right and the left limits $u(t^+):=\lim_{\tau\to
  t^+} u(\tau)$ and $u(t^-):= \lim_{\tau\to   t^-} u(\tau)$ exist. If $u$ is an
energetic solution, we can always modify $u$ such that $u(t)=u(t^+)$ or
$u(t)=u(t^-)$, and we still have an energetic solution. Indeed, for our convex
case we may even set $u(t)=(1{-}\theta)u(t^+)+ \theta u(t^-)$ for any $\theta
\in [0,1]$. Thus, uniqueness holds only if we prescribe the jump behavior,
e.g.\ by asking left continuity, i.e.\ $u(t)=u(t^-)$ for all $t$.

As a consequence of the transfer between the gradient system, for which the
classical results of Br\'ezis \cite[Thm.\,3.1+3.2]{Brez73OMMS} provide existence
and uniqueness, we obtain the following result for the ERIS.

\begin{theorem}\label{th:ExiUni} The ERIS $(X,\calE,\|\cdot\|)$ with
  $\calE(t,u)=t\calJ(u)$ 
  possesses for all initial values $u^0\in X$ with $\calJ(u^0)<\infty$ a unique
  left-continuous energetic solution $u:{[0,\infty[} \to X$ in the sense of
   \eqref{eq:ES.tJ}.
\end{theorem}

For cases with $\calJ(u^0)=\infty$ we still have $\calE(0,u^0)=0$, but 
there is a delicate issue about attainment of the initial condition discussed
in Proposition \ref{pr:Contin.t0} and Remark \ref{rm:Attainment}.\medskip

The plan of the paper is as follows: In Section \ref{se:GS} we recall the
relevant, classical facts from the gradient-flow theory developed in
\cite[Thm.\,3.1+3.2]{Brez73OMMS}. Moreover, we discuss the case that $s \to
\|w'(s)\|$ has a plateau and show that this implies that the solution must
follow a straight line. 

In Section \ref{se:Examples} we consider several examples, first a few simple
finite-dimensional ones. Then, we provide an infinite-dimensional example in 
$\rmL^2(\R)$ where all solutions can be calculated explicitly and where we are
able to choose an initial value $u^0$ such that $\int_{t=0}^1 \|\dot u(t)\| \dd t
=\infty$, i.e.\ the right limit $u(0^+)=\lim_{\tau\to 0^+}u(\tau)$ does not
exist. Finally, we shortly refer to the so-called ``total-variation flow''
that is the main motivation for the study of gradient systems with
one-homogeneous energy. Motivated by questions in image denoising, one considers 
$X= \rmL^2(\Omega)$ and $\calJ(w)=\int_\Omega |\nabla w| \dd x$, see
\cite{BeCaNo02TVFR, AnCaMa04PQEM, BonFig12TVFS}. Also there, solutions with
constant velocity $w'(s)=v_*$ play an important role, e.g.\ in the form 
$w(t,x)= \max\{ 1{-}\lambda_* t, 0 \} \bm 1_\omega$, where $\omega$ is a suitable
subset of $\Omega$, see \cite{BeCaNo02TVFR}. 

In Section \ref{se:GS2ERIS} we discuss how solutions $w$ of the gradient-flow
equation \eqref{eq:GFE} generate energetic solutions $u$, whereas Section
\ref{se:RIS2GF} provides the opposite direction and concludes with the proof
of Theorem \ref{th:ExiUni}.  


\section{The gradient flow}
\label{se:GS}

The theory developed in \cite[Thm.\,3.1+3.2]{Brez73OMMS} can be applied 
to the gradient-flow equation 
\begin{equation}
  \label{eq:GF-J}
  0\in E w'(s)+ \pl\calJ(w(s)), \qquad w(0)=w_0\in X
\end{equation} 
induced by the GS $(X,\calJ,\frac12\|\cdot\|^2)$. 
Since $\calJ$ is non-negative,
 lower semicontinuous, convex, and
has a dense domain, the induced semiflow $(\mfT_s)_{s\geq 0}$ is a
strongly continuous contraction semigroup on all of $X$, i.e.\ 
\begin{equation}
  \label{eq:mfT}
  \begin{aligned}
& \mfT_r \circ \mfT_s = \mfT_{r+s} \text{ for } r,s\geq 0;\\
& s \mapsto w(s)=\mfT_s (u^0) \text{ is a strongly
  continuous solution of \eqref{eq:GF-J} with } w(0)=u^0,\\
&\| \mfT_s(w_1)-\mfT_s(w_0)\| \leq \| w_1{-}w_0\| \quad \text{ for all } w_0,w_1\in X.
\end{aligned} 
\end{equation}
In particular, we have existence and uniqueness for the initial value problem
\eqref{eq:GF-J}.  Every solution satisfies the energy-dissipation balance
\begin{equation}
  \label{eq:EDB-GF}
  \calJ(w(s_2)) + \int_{s_1}^{s_2} \| w'(s)\|^2 \dd s = \calJ(w(s_1))
\text{ for } 0<s_1< s_2, 
\end{equation}
such that $w \in \rmH^1({[s_1,\infty[};X)$ for all $s_1>0$. In
particular, $w'(s)$ exists for a.a.\ $s \geq 0$. 
The last relation in \eqref{eq:mfT} 
applied to $w_0=w(0)$ and $w_1=w(h)=\mfT_h(w(0))$
for $h>0$,  shows that every solution $w(s)=\mfT_s(w(0))$ satisfies 
\[
\| w(s_2{+}h) - w(s_2)\| \leq \| w(s_1{+}h) - w(s_1)\| \text{ for all }
h>0 \text{ and } 0\leq s_1<s_2. 
\]
Dividing by $h$ and taking $h\to 0^+$, we find
$ \| w'(s_2)\| \leq \| w'(s_1)\|$ for a.a.\ $0<s_1<s_2$.  According to
\cite[Thm.\,3.1\,(5)+(6)]{Brez73OMMS}, the one-sided derivative from the right
behaves even better:
\begin{subequations}
  \label{eq:Brez.w'+}
\begin{align}
  \label{eq:Brez.w'+A}
&w'_+(s):= \lim_{h\to 0^+} \frac1h \big(w(s{+}h) - w(s)\big)  
 \quad \text{exists for all } s> 0.
\\ 
  \label{eq:Brez.w'+B}
&s \mapsto w'_+(s) \text{ is continuous from the right}, 
\\
 \label{eq:Brez.w'+C}
& s\mapsto \|w'_+(s)\| \text{ is non-increasing and continuous from the right}. 
\end{align}
\end{subequations}
Since $w'(s)$  exists a.e., we have $w'(s)=w'_+(s)$ for almost all $s> 0$. 

Moreover, it is shown in \cite{AmGiSa05GFMS} that \eqref{eq:GF-J} is
equivalent to the ``Evolutionary Variational Inequality'' (EVI), which here
takes the form 
\begin{equation}
  \label{eq:EVI}
  \frac12\|w(s){-}v\|^2 - \frac12\| w(r){-}v\|^2 \leq (s{-}r)
\big( \calJ(v) - \calJ(w(s))\big) \text{ for all }v\in X \text{ and }
0\leq r< s. 
\end{equation}
Setting $r=0$ and $v=0$ and using $\calJ(0)=0$ we find the a priori estimate  
\[
\calJ(w(s)) \leq \frac1{2s}\| w(0)\|^2 \text{ for } s>0.
\]
Inserting this into the energy-dissipation estimate and using
$\calJ\geq 0$ we find
\[
\frac1{s} \| w(0)\|^2 \geq \calJ(w(s/2))\geq \int_{s/2}^{s}
\|w'(\sigma)\|^2 \dd \sigma \geq \frac s2 \|w_+'(s)\|^2 . 
\]
Thus we conclude the a priori estimate 
\begin{equation}
  \label{eq:Veloc}
   \|w'_+(s)\|  \leq \frac{\sqrt{2}}s \|w(0)\| \quad \text{for all } s>0. 
\end{equation}

So far, we have only used the convexity of $\calJ$ and $\calJ(w)\geq
\calJ(0)=0$. The coming results rely on the assumption that $\calJ$ is
positively 1-homogeneous. 

The first result is relevant because the graph of the 1-homogeneous functional
$\calJ$ contains segments like the rays $\set{\alpha v}{
  \alpha\geq0}$. However, there may be even more segments if the sublevel
$\set{w \in X}{ \calJ(w)\leq 1}$ is not strictly convex. The result
characterizes the case that speed mapping $s\mapsto \|w'_+(s)\|$ is constant on
an interval ${[s_1,s_2[}$. Using the strict convexity of the Hilbert-space norm
$\|\cdot\|$, we conclude that $w$ restricted to the interval $[s_1,s_2]$ must be
affine. This property will be crucial for relating the solutions of the GS
$(X,\calJ,\|\cdot\|)$ to the solutions of the ERIS.

\begin{proposition}[Intervals of constant speed] 
Assume that $w:{[0,\infty[}\to X$ is a solution of \eqref{eq:GF-J}
and that $\| w'_+(s)\| = \mathop{\mafo{const}}$ for $ s \in {]s_1,s_2[}$.
Then,  $\ds w(s)=\frac{s_2{-}s}{s_2{-}s_1} w(s_1)+
\frac{s{-}s_1}{s_2{-}s_1} w(s_2)$ for all $s\in [s_1,s_2]$, i.e.\
$w|_{[s_1,s_2]}$ is affine. 
\end{proposition}
\begin{proof} We let $\phi=\|w'_+(s)\|$ be the constant. Obviously, only
  the case $\phi>0$ is interesting. 

We consider $r,s$ with $s_1 \leq r < s \leq s_2$.
On the one hand, the energy-dissipation balance \eqref{eq:EDB-GF} gives 
\begin{equation}
  \label{eq:J.affine}
  \calJ(w(s)) + ( s{-}r ) \phi^2 = \calJ(w(r)).
\end{equation}
One the other hand, the equation gives $-E w'(s)=\eta(s) \in
\pl\calJ(w(s))$ for a.a.\ $s\geq 0$. Hence, we have $\|\eta(r)\|=\phi$
for a.a.\ $r\in [s_1,s_2]$. Thus, convexity of $\calJ$ gives the lower
estimate 
\[
\calJ(w(s)) \geq \calJ(w(r)) + \langle \eta(r),w(s){-}w(r)\rangle 
\geq \calJ(w(r)) -  \|\eta(r)\| \|w(s){-}w(r)\|. 
\] 

Together with \eqref{eq:J.affine}, this implies 
\[
\phi \|w(s){-}w(r)\| \geq \calJ(w(r)) - \calJ(w(s)) = (s{-}r) \phi^2. 
\]
Combining this with the trivial upper bound 
\[
\|w(s)- w(r)\| \leq \int_r^s \|w'_+(\sigma)\| \dd \sigma \leq (s{-}r) \phi,
\]
we conclude $\| w(s){-}w(r)\| = (s{-}r) \phi $, which implies that
$w_{[s_1,s_2]}$ is a geodesic curve in the Hilbert space
$(X,\|\cdot\|)$, which implies that it is a straight line. 
\end{proof}

We continue with some auxiliary result for the solutions
of the gradient system that will be needed later.

\begin{proposition}[Norm and energy decay]\label{pr:Norm.J.decay}
Assume that $w:{[0,\infty[} \to X$ is a solution of \eqref{eq:GFE}. Then,
$s\mapsto \| w(s)\|$ is non-increasing. More precisely, we have 
\begin{equation}
  \label{eq:wNorm.decay}
 \frac12 \|w(r)\|^2 + \int_s^r \calJ(w(\sigma))\dd
  \sigma =  \frac12\|w(s)\|^2  \text{ for } 0 < s <r. 
\end{equation}
In particular, the energy $s\mapsto \calJ(w(s))$ is integrable  with 
\begin{equation}
  \label{eq:J.integrable}
  \int_0^\infty \calJ(w(s)) \dd s \leq \frac12 \| w(0)\|^2,
\end{equation}
which is non-trivial for $s\approx 0$ as well as for $s\to \infty$.
\end{proposition}
\begin{proof}
Since $\calJ$ is 1-homogeneous, we have $\langle \eta, w\rangle
=\calJ(w)$ for all $\eta \in \pl\calJ(w)$.  Hence, we may multiply
\eqref{eq:GFE} by $w(s)$ to obtain  
\[
\frac{\rmd}{\rmd s} \frac12\|w(s)\|^2 = \langle w(s),w'(s)\rangle = -
\calJ(w(s)) \leq 0. 
\]
Thus, for every solution $w$ of \eqref{eq:GFE} and  $0<s_1<s_2$ we have 
\[
\int_{s_1}^{s_2} \calJ(w(s)) \dd s = -\int_{s_1}^{s_2} \frac{\rmd}{\rmd
s} \frac12\|w(s)\|^2 \dd s = \frac12 \|w(s_1)\|^2 - \frac12
\|w(s_2)\|^2,
\]
which is the desired result.
\end{proof}

The previous result can be seen as a special case of the ``inverse
energy-dissipation balance'' also used in \cite[Lem.\,4.1]{HeNeVa19?SHLC}. For
quadratic dissipation potentials $\calR(v)=\frac12\|v\|^2$ in a Hilbert space
and general convex functionals $\calJ$, the gradient-flow equation
$0=\rmD\calR(w')+\pl\calJ(w)$ can also rewritten as
\[
\calR(w(s)) + \int_r^s \big\{\calJ(w(\sigma)) +\calJ^*({-}\rmD\calR(w'(\sigma)
) \big\} \dd \sigma = \calR(w(r)) \text{ for all } s>r\geq 0. 
\]
In our special case of the 1-homogeneous $\calJ$ the Legendre-Fenchel dual
$\calJ^*$ satisfies $\calJ^*(\xi)=0$ for $\xi \in \pl\calJ(0)$ and
$\calJ^*(\xi)=\infty$ otherwise . Hence, $\calJ^*$ doesn't show up in
\eqref{eq:wNorm.decay}.

We conclude this section by an estimate on the extinction
time $S_\mafo{extinct}(u^0)$, i.e.\ the solution $w$ with $w(0)=u^0$ satisfies
$w(s)=0 $ for $s\geq S_\mafo{extinct}(u^0)$. 

\begin{proposition}[Extinction] Assume that for the GS
  $(X,\calJ,\frac12\|\cdot\|^2)$  we additionally have
  \begin{equation}
    \label{eq:ExtinctAssump}
    \calJ \text{ is 1-homogeneous and } \quad \exists\,\beta>0\ \forall \, u\in
    X: \ \calJ(u) \geq \beta\|u\|.
  \end{equation}
Then, the solution $w$ with $w(0)=u^0$ satisfies $\|w(s)\|\leq \max\big\{ 0, 
 \|u^0\|{-}\beta s\big\}$, which implies $S_\mafo{extinct}(u^0)\leq \|u^0\|/\beta$.  
\end{proposition} 
\begin{proof} We simply estimate the norm via $\frac12\frac\rmd{\rmd s}
  \|w(s)\|^2 = \langle Ew',w\rangle =- \calJ(w) \leq -\beta \|w(s)\|$. Setting
  $\rho(s)=\|w(s)\|\geq 0$, this means $\rho(\rho'{-}\beta)\leq 0$, which
  implies the desired result.  
\end{proof}

\section{Some examples}
\label{se:Examples}

We first consider three examples with $X=\R^2$ equipped with the Euclidean
norm and then a simple case on $X=\rmL^2(\Omega)$ that can be solved
explicitly. 

\subsection{Piecewise affine example} 
\label{su:PieceAff}

We consider $\calJ(u)=\max\{|u_1|,2|u_2|\}$. For the subdifferential of $\calJ$
we find
\[
\pl\calJ(u)= \left\{ \ba{cl}
 \{ (\sign(u_1),0)^\top\}&\text{for }|u_1|>2|u_2|,\\[0.2em]  
 \{ (0,2\sign(u_2))^\top\}&\text{for }|u_1|<2|u_2|,\\[0.3em] 
\bigset{ (\sign(u_1)(1{-}\theta),2\sign(u_2)\theta )^\top}{\theta \in
  [0,1]}& \text{for }|u_1|=2|u_2|\neq 0,\\[0.3em]   
 \bigset{(\eta_1,\eta_2)^\top}{2|\eta_1|{+}|\eta_2|=2}& \text{for }u=0.  \ea \right.
\]
For the case $|u_1|=2|u_2|\neq 0$, the minimal element $\pl^0\calJ(u)$ of
$\pl\calJ(u)$ takes the form
$ \pl^0\calJ(u)= \frac15 (4\sign(u_1),2\sign(u_2))^\top$, and all other cases
are trivial.

\begin{figure}
\centerline{
\begin{minipage}{0.45\textwidth}
\begin{tikzpicture}
\draw[thick,->] (-3.2,0)--(3.2,0) node[right]{$u_1$};
\draw[thick,->] (0,-2)--(0,2) node[right]{$u_2$};
\foreach \mm in {0.3,0.6,0.9,1.2,1.5}
 \draw[gray!50, very thick] ({-2*\mm},{-\mm}) rectangle  ({2*\mm},\mm);
\draw[blue, ultra thick, ->] (-3.0,0)--(-0.9,0);
\draw[blue, ultra thick] (-1.2,0)--(0,0);
\draw[blue, ultra thick,] (3.0,1)--(2,1)--(0,0);
\draw[blue, ultra thick,>->] (2.5,1)--(2,1)--(1,0.5);
\draw[blue, ultra thick,] (1.6,-1.5)--(1.6,-0.8)--(0,0);
\draw[blue, ultra thick,>->] (1.6,-1.2)--(1.6,-0.8)--(0.8,-0.4);
\draw[blue, ultra thick,] (-3,1.5)--(0,0);
\draw[blue, ultra thick,->] (-3,1.5)--(-1.6,0.8);
\draw[blue, ultra thick,] (0,-1.5)--(0,0);
\draw[blue, ultra thick,->] (0,-1.5)--(0,-0.8);
\end{tikzpicture}
\end{minipage}
\quad
\begin{minipage}{0.35\textwidth}
\caption{Gradient flow in $\R^2$ for $\calJ(u)=\max\{|u_1|,2|u_2|\}$: the level
  sets of $\calJ$ are indicated in gray, and a few orbits are drawn in blue.}
\label{fig:Affine}
\end{minipage}}
\end{figure}
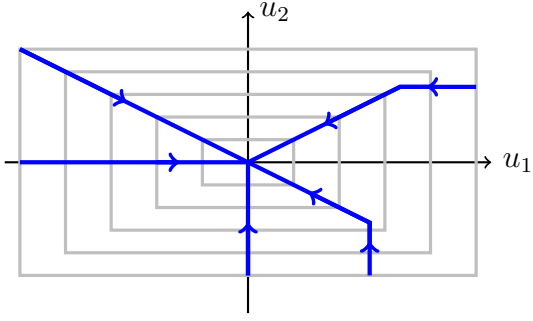

Since $\pl^0\calJ$ only takes finitely many values, the solutions are easily
constructed by straight lines in $\R^2$. Without loss of generality we
consider the case $w(0)=u^0=(u_1^0,u_2^0)$ with $0<2u_2^0< u_1^0$, 
then the explicit solution reads 
\[
w(s) = \left\{ \ba{cl} u^0 - (s,0)^\top& 
 \text{for } s \in [0,u_1^0{-} 2u_2^0], \\[0.3em] 
\frac{2u_1^0{+}u_2^0{-}2s}{5} \:(2,1)^\top& 
  \text{for }s \in [u_1^0{-} 2u_2^0,u_1^0{+}u_2^0/2], 
\\[0.3em] 0 & \text{for }s\geq u_1^0{+}u_2^0/2. \ea \right.
\] 
We see that $w'_+$ only takes three values, namely $w'_+(s)\in \{(1,0)^\top,
(4/5,2/5)^\top,(0,0)^\top\}$. Moreover, every solution reaches
$w=0$ in finite 
time, namely $S_\text{extinct}(u^0)=|u_1^0|+|u_2^0|/2$, see also Figure
\ref{fig:Affine}.  

\subsection{Singular potential}
\label{su:ExaSing}

We again consider the Hilbert space $\R^2$ with the Euclidean
norm. The potential is 
\[
\calJ(u)= \frac{|u_1|^{\alpha+1}}{u_2^\alpha} \text{ for } u_2>0, \quad
\calJ(0)=0,\quad\text{ and }\calJ(u)=\infty \text{ otherwise},
\]
where the exponent $\alpha$ satisfies $\alpha\geq 1$, such that
$\calJ$ is indeed convex and lower semicontinuous. The space $X$ is
now the closure of the domain of $\calJ$, namely $X=\R\ti {[0,\infty[} \subset \R^2$. 

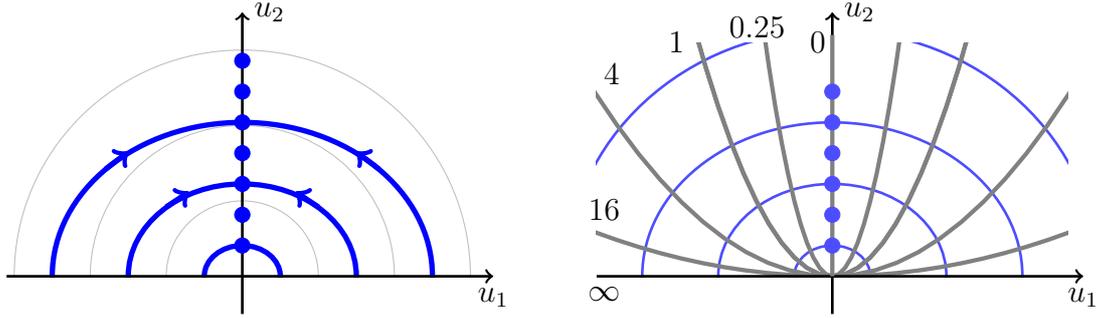
\begin{figure}[h]
\usetikzlibrary{calc}
\mbox{}\hfill
\begin{tikzpicture}
\begin{scope}
\clip (-3.1, -0.01) -- (-3.1,3.1) --(3.1,3.1) --(3.1,-0.01)--cycle;
\draw[color=gray!50 ] (0,0) circle (1);
 \draw[color=gray!50 ] (0,0) circle (2);
\draw[color=gray!50 ] (0,0) circle (3);
\draw[line width=2pt,color=blue] (0, 0) ellipse (0.5 and 0.408);
\draw[line width=2pt,color=blue] (0, 0) ellipse (1.5 and 1.225);
\draw[line width=2pt,color=blue] (0, 0) ellipse (2.5 and 2.041);
\end{scope} 
\draw[thick,->,line width=1pt] (0,-0.5) -- (0,3.5)node[right]{$u_2$};
\draw[->,line width=1pt] (-3.1,0) --(3.3,0)node[below]{$u_1$};

\draw[color=blue, fill ] (0.0,0.408) circle (0.1); 
\draw[color=blue, fill ] (0.0,0.816) circle (0.1); 
\draw[color=blue, fill ] (0.0,1.225) circle (0.1); 
\draw[color=blue, fill ] (0.0,1.633) circle (0.1); 
\draw[color=blue, fill ] (0.0,2.041) circle (0.1); 
\draw[color=blue, fill ] (0.0,2.449) circle (0.1); 
\draw[color=blue, fill ] (0.0,2.858) circle (0.1); 

\draw[->, color=blue, line width=2pt] (-1.8,1.42) -- (-1.5,1.65); 
\draw[->, color=blue, line width=2pt] (1.8,1.42) -- (1.5,1.65); 
\draw[->, color=blue, line width=2pt] (1,0.92) -- (0.7,1.12); 
\draw[->, color=blue, line width=2pt] (-1,0.92) -- (-0.7,1.12); 
\end{tikzpicture} \hfill
\begin{tikzpicture}
\begin{scope}
\clip (-3.1, -0.01) -- (-3.1,3.1) --(3.1,3.1) --(3.1,-0.01)--cycle;
\draw[line width=1pt,color=blue!70] (0, 0) ellipse (0.5 and 0.408);
\draw[line width=1pt,color=blue!70] (0, 0) ellipse (1.5 and 1.225);
\draw[line width=1pt,color=blue!70] (0, 0) ellipse (2.5 and 2.041);
\draw[line width=1pt,color=blue!70] (0, 0) ellipse (3.5 and 3.265);
\draw[gray, ultra thick,domain = -2:2 ] plot (\x, {pow(\x,2)});
\draw[gray, ultra thick,domain = -1:1 ] plot (\x, {4.0* pow(\x,2)});
\draw[gray, ultra thick,domain = -2:2 ] plot (\x, {pow(\x,2)});
\draw[gray, ultra thick,domain = -3.4:3.4] plot (\x,{0.25*pow(\x,2)});
\draw[gray, ultra thick,domain = -3.4:3.4] plot (\x,{0.06*pow(\x,2)});
\end{scope} 
\draw[thick,->,line width=1pt] (0,-0.5) -- (0,3.5)node[right]{$u_2$};
\draw[->,line width=1pt] (-3.1,0) --(3.3,0)node[below]{$u_1$};
\draw[gray, ultra thick] (0,0)-- (0,3.2);
\node[below] at (-3,-0.0) {$\infty$};
\node[above] at (-3,0.6) {$16$};
\node[above] at (-2.9,2.4) {$4$};
\node[left] at (-1.8,3.1) {$1$};
\node at (-.99,3.3) {$0.25$};
\node at (-0.19,3.1) {$0$};
\draw[color=blue!70, fill ] (0.0,0.408) circle (0.1); 
\draw[color=blue!70, fill ] (0.0,0.816) circle (0.1); 
\draw[color=blue!70, fill ] (0.0,1.225) circle (0.1); 
\draw[color=blue!70, fill ] (0.0,1.633) circle (0.1); 
\draw[color=blue!70, fill ] (0.0,2.041) circle (0.1); 
\draw[color=blue!70, fill ] (0.0,2.449) circle (0.1); 
\end{tikzpicture}
\hfill\mbox{}
\caption{Solutions (right) and level sets of potential $\calJ$ (left) 
for $\alpha=2$.}\label{fig:SingPot} 
\end{figure}

The point is that we are able to characterize the solutions
explicitly. Using the gradient-flow equation 
\[
w'_1 = - (\alpha{+}1) \frac{|w_1|^{\alpha-1}w_1}{w_2^{\alpha}}
, \qquad 
w'_2= \alpha \frac{|w_1|^{\alpha+1}}{w_2^{\alpha + 1} }
\]
we easily see that the function 
\[
\Phi_\alpha(w)= \frac{\alpha}{\alpha+1} \, w_1^2 +  w_2^2 
\]
is constant along solutions. All solutions have $w_1 w'_1\geq0$ and
$w'_2 \geq 0$, see Figure \ref{fig:SingPot}.

The first observation concerns the limit $\lim_{s\to \infty} w(s)=
\lim_{t\to \infty}u(t)$ for the unique solution starting in
$u^0=(u_1^0,u_2^0)$ with $u_2^0\geq 0$. Using the first integral
$\Phi$ and the signs of $w'_j$, we find that the the solution converges
to limit point 
\[
\mathfrak L(u^0):= \lim_{s\to \infty} w(s)=
\lim_{t\to \infty}u(t) \ = \ \Big(\, 0\, , \sqrt{\Phi_\alpha(u^0)}\Big). 
\]
As expected, the mapping $u^0\mapsto \mathfrak
L(u^0)$ is 1-homogeneous, but otherwise it is nonlinear. 

We may also analyze the decay properties of $w(s)$ for $s\to \infty$
or of $u(t)$ for $t\to \infty$. Restricting to the level set
$\Phi(w(s)) = \Phi(u^0)$, it is sufficient to solve an ODE for $w_1$ of
the form $w'_1(s)=F_{u^0}(w_1)$, where $F_{u^0}(w_1)=-
c_{u^0}|w_1|^{\alpha -1} w_1 +O(|w_1|^{\alpha+1})$. Hence, in the case
$\alpha>1$ we find algebraic decay of the form 
\[
w(s)- \mathfrak L(u^0) \approx c
s^{-1/(\alpha-1)} \text{ for } s \to \infty. 
\]
Using $t(s)=\|w'(s)\|$, we hence find $u(t)-\mathfrak L(u^0)\approx c
t^{-1/\alpha}$.   

The case $\alpha=1$ is special, since $w(s)$ convergence exponentially
to $\mathfrak L(u^0)$, while $u(t) - \mathfrak L(u^0) \approx c
/t$.\medskip

It is also interesting to analyze the behavior of solutions starting
with $\calJ(u^0)=\infty$, i.e.\ $u^0=(\beta,0)$ with $\beta\neq 0$. 
Using the first integral $\Phi$ again, we can now write an ODE for $w_2$,
namely 
\[
w'_2 = G_{u^0}(w_2) = c_{\alpha,u^0}w_2^{-(\alpha+1)} +
O(w_2^{-\alpha}) \text{ for } w_2\to 0.
\]
Thus, we find $w(s)- u^0 \approx c s^{1/(\alpha+2)}$  and thus 
\[
\| w'(s)\| \approx s^{-(\alpha+1)/(\alpha+2)}\quad \text{ and } \quad 
\calJ(w(s)) \approx s^{-\alpha/(\alpha+2)}\quad  \text{ for }s\to 0. 
\]
Note that $\|w'(\cdot)\|$ and $\calJ(w(\cdot))$ are integrable near
$s=0$, while $s\mapsto \| w'(s)\|^2$ is not. 

For the energetic solution we find using $t(s)=1/\|w'(s)\| \approx
s^{(\alpha+1)/(\alpha+2)} $ the relations 
\[
u(t)-u^0\approx  t^{1/(\alpha+1)},\quad \| \dot u(t)\| \approx
t^{-\alpha/(\alpha+1)} , \quad \calJ(u(t))\approx
t^{-\alpha/(\alpha+1)}.
\]
In particular, we conclude that $t\mapsto
\calE(t,u(t))=t\,\calJ(u(t))$ is continuous and $\|\dot u\|$ and
$\calJ(u(\cdot))$ are integrable.

\subsection{Degenerate potential}
\label{su:Degenerate}

The functional $\calJ(u)=u_1+\|u\|=\sqrt{u_1^2{+}u_2^2}-u_1$ has the
property that $\calJ(u)=0$ if and only if $u\in \set{(a,0)}{a\geq 0}$.   
For $w(0)=(a,0)$ we have $w(s)=(a,0)$ if $a\geq 0$ and
$w(s)=(\min\{a{+}2s,0\}, 0)$ for $a <0$. 

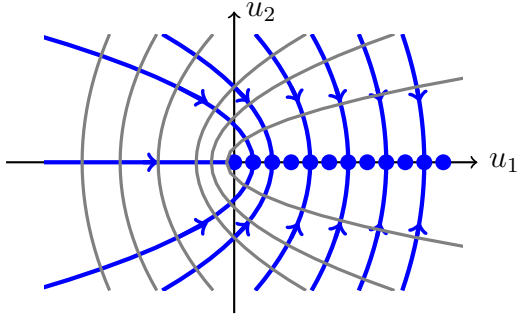
\begin{figure}
\centerline{
\begin{minipage}{0.5\textwidth}
\begin{tikzpicture}
\draw[thick,->] (-3,0)--(3.2,0) node[right]{$u_1$};
\draw[thick,->] (0,-2)--(0,2) node[right]{$u_2$};
\begin{scope}
\clip (-2.5,-1.7) rectangle (3.0,1.7);
\draw[blue, ultra thick, ->] (-3,0)--(-1,0);
\draw[blue, ultra thick] (-1.2,0)--(0,0);
\foreach \i in {0,...,11} 
  \fill[blue] (\i/4,0) circle (3pt);
\foreach \mm in {0.5,1,2,3,4,5}
 \draw[blue, ultra thick, domain =-1.7:1.7] plot ({(\mm^2-\x*\x)/(2*\mm)},\x);
\foreach \mm in {0.5,1,2,3,4,5}
 \draw[blue, ultra thick,->] ({(\mm^2-1)/(2*\mm)},1.0)--
                             ({(\mm^2-1)/(2*\mm)+0.2/\mm}, 0.8); 
\foreach \mm in {0.5,1,2,3,4,5}
 \draw[blue, ultra thick,->] ({(\mm^2-1)/(2*\mm)},-1.0)--
                             ({(\mm^2-1)/(2*\mm)+0.2/\mm}, -0.8); 
%
\foreach \mm in {0.2,0.6,1,2,3,4}
 \draw[gray, very thick, domain =-1.7:1.7] plot ({(\x*\x-\mm*\mm)/(2*\mm)},\x);
\end{scope}
\end{tikzpicture}
\end{minipage}
\begin{minipage}{0.3\textwidth}
\caption{Gradient flow in $\R^2$ for $\calJ(u)=\|u\|-u_1$. Level sets of
  $\calJ(\cdot)=c \in \{0.1,0.3,0.5,1,1.5,2\}$ are shown
  in gray, and the orbits are drawn in blue.}
\label{fig:Degenerate}
\end{minipage}
}
\end{figure}

Moreover, it is easy to see that $M(u)=\sqrt{u_1^2{+}u_2^2}+u_1$ is conserved
along solutions. Hence, the orbits lie on the level sets $M(u)=m\geq0$, which
are parabolas, and satisfy $u_1(t)=(m^2{-}u_2(t)^2)/(2m)$, see Figure 
\ref{fig:Degenerate}. 
In particular, a solution starting in $u^0 \in \R^2$ satisfies
$u(t)\to \big(\frac12M(u^0),0\big)$ for $t \to \infty$, where the convergence
is exponential for $u^0_2\neq 0$.

\subsection{An infinite dimensional example}
\label{su:InfiniteExa}

We consider the Hilbert space $X= \rmL^2(\R)$ and the 1-homogeneous functional
$\calJ(u)= \int_\R a(x) |u(x)| \dd x $ for some non-negative and measurable
function $a: \R \to {[0,\infty[}$. We can give the subdifferential $\pl\calJ$
and the minimal elements $\pl^0\calJ(u)$ explicitly via
\begin{align}
 \label{eq:infExa.plJ}
 \pl\calJ(u)&
 = \bigset{\xi \in  \rmL^2(\R) }{ \xi(x) \in a(x)\Sign(u(x)) \text{ a.e.\;}}
   \text{ and }  \pl^0\calJ(u)= a(x)\sign(u(x)).
\end{align}

The associated gradient-flow equation reads 
\[
0 \in w'(s,x) + a(x) \mathop{ \mafo{Sign}}(w(s,x)), \qquad w(0,x)=u^0(x).
\]
We easily see that the solutions are given by the explicit formula
\[
w(s,x) =\mathop{\mafo{Sign}}(u^0(x))\:\max\big\{ \,|u^0(x)| - a(x)s\,, \;0 \, \big\}
\]
for the solutions $w(s)$.  For the time derivative we obtain the formula
\[
w'(s,x)=- a(x) \mathop{\mafo{Sign}}(u^0(x))\:\bm1_{\{x : 
|u^0(x)| - a(x)s>0\}} (x)
\]
and see that $w'(s)$ can be constant on an interval ${]s_1,s_2[}$ if the 
image of the function $x \mapsto |u^0(x)|/a(x)$ intersected with ${]s_1,s_2[}$
has Lebesgue measure $0$. 

To see one typical behavior for $s\approx 0$ and $s\gg1$ we look at the simple
case 
\[
a \equiv 1 \qquad \text{and} \qquad u^0(x) =\min\{|x|^{-\alpha}, |x|^{-\beta}\}
\ \text{ with } 0<\alpha < 1/2< \beta <1. 
\]
Since $u^0$ is even and strictly decreasing for $x>0$, it is easily seen that
the support of $w(s,\cdot)$ is given by $[-X(s),X(s)]$ with
$X(s)= \min\{ s^{-1/\alpha}, s^{-1/\beta}\}$ and
\[
w(s,x)= \big(u^0(x)-u^0(X(s))\big) \bm1_{\{|x|\leq X(s)\}}(x) .
\]
For $s\to 0$ we obtain $\calJ(w(s))=\frac{2\beta}{1-\beta} s^{-(1-\beta)/\beta} +
O(1)$, while for $s\geq 1$ we have $J(w(s))=\frac{2\alpha}{1-\alpha}
s^{-(1-\alpha)/\alpha} $, which is compatible with $\calJ(w(s))\dd s <\infty$,
cf.\ \eqref{eq:J.integrable}. 

For the velocity and the slope we have
$\|w'(s)\|= \|\pl^0\calJ(w(s))\|_*= (2X(s))^{1/2}$, which is compatible with
$\|w'(s)\|\leq C/s$, cf.\ \eqref{eq:Veloc}.  Moreover, by choosing $u^0$
suitably, we can find a support mapping $s\mapsto X(s)$ such that one or both
of the integrals $\int_0^1\|w'(s)\|\dd s$ and $\int_1^\infty\|w'(s)\|\dd s$ are
infinite.
 
\subsection{Total variation flow}
\label{su:TotalVariation}

An important motivation for the present work is the so-called total variation
flow, which plays an important role in image processing. We refer to
\cite{RuOsFa92NTVB, BeCaNo02TVFR, BonFig12TVFS, KiMuRy13ACST, BoDuMa15TDVA}
or the monograph \cite{AnCaMa04PQEM} and the references therein. 

According to \cite{RuOsFa92NTVB} the denoising of an image $f$ given over a domain
$\Omega \subset \R^2$ can be done by considering the $\rmL^2$ gradient flow for
the convex functional 
\[
\wt \calJ: u \mapsto \int_\Omega \big\{ |\nabla u| +  
 \frac\kappa p | u{-}f|^p \big\} \dd x\,.
\] 
This leads to the parabolic equation 
\[
w'= \DIV\Big( \frac1{|\nabla u|}\,\nabla u \Big) - \kappa |u {-}f|^{p-2} (u{-}f)
\qquad \text{ in } \Omega,
\]
which has to be completed by no-flux boundary conditions and 
interpreted in a suitable weak sense. 

Obviously, our theory only applies in the case $f\equiv 0$ and $\kappa =0$ or
$\kappa>0$ and $p=1$. The former case is also relevant in crystal growth, see
\cite[Eqn.\,(61)]{QuaMar08ADCR} and \cite{KiMuRy13ACST}.  The latter work
addresses in particular the one-dimensional case and shows that facets are
preserved.  As a consequence for $\Omega=\R$, the class of step functions is
invariant under the gradient flow. Choosing points $y_0<y_1< \ldots < y_N$ and
setting 
\[
w(s,x) = \sum_{i=1}^N \alpha_i(s) \bm1_{{]y_{i-1},y_i]}}(x) ,
\]
we can reduce the $\rmL^2$ norm and the total variation functional to obtain
\begin{align*}
\|\bfalpha(s)\|_\bfy&:= 
\|w(s)\|_{\rmL^2} = \Big( \sum_{i=1}^N \alpha_i(s)^2 (y_i{-}y_{i-1})\Big)^{1/2},
\\
\bfJ(\bfalpha(s))&:=\calJ(w(s))= \int_\R |\pl_x w(s,x)| \dd x =|\alpha_1(s)| +
|\alpha_N(s)| +  \sum_{i=2}^{N} \big| \alpha_i - \alpha_{i-1}\big|.
\end{align*}
Thus, the evolution of the vector $\bfalpha(s)=(\alpha_i(s))_i \in \R^N$ is
indeed determined by the finite-dimensional gradient-flow equation for the
gradient system $(\R^N,\bfJ,\|\cdot\|_\bfy)$. It takes the form
\begin{equation}
  \label{eq:TV.FiniDim}
  (y_i{-}y_{i-1})\,\alpha_i'(s) \in 
          - \Sign\big(  \alpha_i(s){-}\alpha_{i-1}(s)\big)
          - \Sign\big(  \alpha_i(s){-}\alpha_{i+1}(s)\big)  \ \text{ for }i=1,...,N,
\end{equation}
where we set $\alpha_0(s)=\alpha_{N+1}(s)=0$ and use the set-valued function
$\Sign$ with $\Sign(0)=[-1,1]$. We refer to  \cite{BonFig12TVFS,
  KiMuRy13ACST} for illustrative examples.

\section{Energetic solutions} 
\label{se:ConstES}

Energetic solutions are defined by stability and by the energy balance in
\eqref{eq:def.ES}. We refer to \cite{Miel05ERIS} and \cite{MieRou15RIST} for an
introduction and a more extensive theory, respectively. 
For definiteness we rewrite the definition of energetic solutions for our special
ERIS with $\calE(t,u)=t\calJ(u)$. We call a mapping
$u:{[0,\infty[} \to X$ an energetic solution for the ERIS
$(X,\calE,\|\cdot\|)$, 
if $t\mapsto \calJ(u(t))$ lies in $\rmL^1_\mafo{loc}({[0,\infty[})$ and if for 
all $r,t\geq 0$ with $r<t$ we have
\begin{subequations}
 \label{eq:ES.tJ}
 \begin{align}
  \text{(S)}\quad &\label{eq:ES.tJ.S} 
  \forall\, \wt u\in X: \quad t\calJ(u(t)) \leq t\calJ(\wt u)+ \| \wt u{-}u\|;
 \\
  \text{(E)}\quad &\label{eq:ES.tJ.E}
  t\calJ(u(t)) + \mafo{Var}_{\|\cdot\|} (u;[r,t]) = r\calJ(u(r))
  +\int_r^t \calJ(u(\tau)) \dd \tau;
 \\
 \text{(I)}\quad &\label{eq:ES.tJ.I}
  u(0)=u(0^+):=\lim_{\tau\to 0} u(\tau) \ \text{ in the norm topology}.
 \end{align}
\end{subequations}
Note that we ask the energy balance (E) on all compact subintervals $[r,t]$ 
of ${[0,\infty[}$, whereas it is usual to impose it only on $[0,T]$ for all
$T>0$. However, we have a singular situation at $t=0$, because of 
$\calE(0,u)=0\calJ(u)=0$. 

Here it is important that $\int_0^t \calJ(u(\tau))\dd \tau <\infty$ implies
that $\mafo{Var}_{\|\cdot\|} (u;[0,t])<\infty$ for all $t>0$, see Proposition
\ref{pr:Contin.t0}. This means that
the limit from the right  $u(0^+)=\lim_{\tau\to 0^+} u(\tau)$ exists and the
attainment of the initial condition in (I) is well-defined.

\subsection{Preliminaries on ERIS} 
\label{su:PreERIS}

In general, the stability condition (S) in \eqref{eq:def.ES} is best formulated
via the sets of stable states
\[
\calS(t):=\bigset{u\in X}{ \forall\,v\in X: \ \calE(t,u)\leq
  \calE(t,v)+ \| v{-}u\|  }.
\]
For (S) in \eqref{eq:ES.tJ.S} we use the convexity of $\calJ$ and 
$\calE(t,u)=t\calJ(u)$ to find
\[
\calS(t)=\bigset{u\in X}{ t\pl\calJ(u) \in B_1(0)\subset X^*} = \bigset{u\in
  X}{ \|\pl^0\calJ(u)\|_* \leq 1/t}.   
\]
Hence, for $t>0$ the stability of $u$ implies that $u \in
\mafo{dom}(\pl\calJ)$. We obviously have 
\begin{equation}
  \label{eq:Stab.Sets}
  X=\calS(0) \ \supset \ \calS(t_1)\ \supset \  \calS(t_2) \ \supset \
  \set{u \in X}{ \calJ(u)=0 } \quad  \text{ for } 
0\leq t_1< t_2. 
\end{equation}
Moreover, each $\calS(t)$ is a cone, i.e.\ $\lambda \geq 0$ and $u \in
\calS(t)$ implies $\lambda u \in \calS(t)$. However, in general 
these sets are not convex and not weakly
closed. Indeed for the example $\calJ(u)=\max\{|u_1|,2|u_2|\}$ from Section
\ref{su:PieceAff} we have 
\begin{equation}
  \label{eq:S.notcvx}
  \calS(t)=\left\{ \ba{cl}\R^2& \text{for }t \in [0, 1/2],\\[0.3em]
\bigset{(u_1,u_2)}{ |u_1|\geq 2|u_2| }& \text{for } t \in {]1/2,1]},\\[0.3em]
\bigset{(u_1,u_2)}{ |u_1|= 2|u_2| } & \text{for } t \in {]1,2/\sqrt5]},\\[0.3em]
\{(0,0)\}& \text{for } t > 2/\sqrt5.   \ea \right.
\end{equation}
However, the stability sets are strongly closed, which follows easily from the
lower semicontinuity of $\calJ$, see  also \cite[Prop.\,5.9]{Miel05ERIS}.

\begin{lemma}[Strong closedness of stability sets]
\label{le:S-closed} 
For a sequence $(t_n,u_n)\in {[0,\infty[} \ti X$ we have
\begin{equation}
  \label{eq:S-t-closed}
 \Big( u_n\in \calS(t_n) \ \text{ and }\  (t_n,u_n)\to (t,u)\Big)\quad
  \Longrightarrow \quad u\in \calS(t). 
\end{equation}
\end{lemma} 

The following example, which is the rate-independent analog of 
the GS studied in  Section   \ref{su:InfiniteExa}, 
provides a non-trivial case, in which we are able to show that the case 
$\mafo{Var}_{\|\cdot\|}(u;[0,1])=\infty$ may actually occur, which implies that
the limit $u(0^+)=\lim_{\tau\to 0^+}u(\tau)$ does not exist and the attainment
of the initial condition $u(0)=u(0^+)$ doesn't make sense.

\begin{example}[Case with $\int_0^1\calJ(u(\tau))\dd\tau=\infty$] 
\label{ex:int.J.infty}\slshape 
We return to the example treated in Section
  \ref{su:InfiniteExa} with $X=\rmL^2(\R)$ and $\calJ(u)=\int_\R |u(x)|\dd x$
  (i.e.\ $a\equiv 1$). Starting with a non-negative, even $u^0$, such that
  $u^0|_{{]0,\infty[}}$ is differentiable and strictly decreasing, we obtain the
  energetic solution 
\[
u(t,x)= \max\{ 0,\: u^0(x) - S(t)\} \quad \text{with } S(t)=u^0(X(t)) \text{
  and } X(t)=1/(2t^2).
\]  
To understand the construction, consider the special case
$u^0(x)=(2|x|)^{-\beta}$ for $|x|\geq 1$ with $\beta>1/2$ to have $u^0\in
\rmL^2(\R)$. Then $S(t)=t^{2\beta}$. 

To show that $u$ is an energetic solution we observe that $\pl^0\calJ(u)$ given
in \eqref{eq:infExa.plJ} reads
\[
\pl^0\calJ(u(t))(x)=\bm1_{[-X(t),X(t)]}(x) \quad \text{giving } 
\| \pl^0\calJ(u(t))\|_*= \big(2X(t)\big)^{1/2}= 1/t,
\]
which implies $u(t)\in \calS(t)$ and (S) in \eqref{eq:ES.tJ} is satisfied. 
 
Moreover, $t \mapsto u(t)$ is differentiable with 
\[
\dot u(t,\cdot) = -\dot S(t)
\bm1_{[-X(t),X(t)]} \quad \text{giving }\|\dot u(t)\| 
= |\dot S(t)|\big(2X(t)\big)^{1/2} =  |\dot S(t)|/t .
\] 
Hence, on the one hand the variation can be calculated to obtain
\begin{align*}
\mafo{Var}_{\|\cdot\|}(u;[r,t]) &= \int_r^t \|\dot u(\tau)\| \dd \tau
= \int_r^t \frac{\dot S(t)}t\dd t 
= \frac{S(t)} t - \frac{S(r)}r + \int_r^t \frac{S(\tau)}{\tau^2} \dd \tau.
\end{align*}
On the other hand, the functional $\calJ$ can be evaluated explicitly via
\begin{align*}
\calJ(u(t))&= 2 \int_0^{X(t)}\!\big\{u^0(x)-S(t)\big\} \dd x 
  \ \overset{x=X(\tau)}= \ 2\int_t^\infty \!\big\{ u^0(X(\tau)) - u^0(X(t))\big\}
  |\dot X(\tau)| \dd \tau 
\\
&= 2\int_t^\infty \frac{S(\tau)-S(t)}{\tau^3} \dd \tau = 2\int_t^\infty
\frac{S(\tau)}{\tau^3} \dd \tau  - \frac{S(t)}{t^2}. 
\end{align*}
With this, the energy balance (E) in \eqref{eq:ES.tJ} follows by a
straightforward calculation.  

For this case we can construct an example where $\int_0^1 \calJ(u(t))\dd t
=\infty$ and hence $\mafo{Var}_{\|\cdot\|}(u;[0,1])=\infty$ as well. In
particular, $\lim_{\tau \to 0^+} u(\tau)$ does not exist. To obtain such
a case choose $u^0$ with $u^0(x)= |x|^{-1/2} (\log|x|)^{-\beta}$ for $|x|\geq
x_*\gg1$ and $u^0(x)=0$ otherwise. Then $u^0\in \rmL^2(\R)$ for $\beta >1/2$. 
For $t\to 0$ we obtain $S(t) \approx c_0 t (\log(1/t))^{-\beta}$ and $\dot S(t) =
c_1 (\log(1/t))^{-\beta}$. We conclude 
\[
 \mafo{Var}_{\|\cdot\|}(u;[0,1]) = \int_0^1 \frac{\dot S(t)}t \dd t \approx
 \int_0^t \frac{\dd t}{t\,(\log(1/t))^\beta} = \infty \quad \text{for }\beta \leq
 1.
\] 
Thus, for $\beta \in {]1/2,1]}$ we obtain a case where $\lim_{\tau \to
  0^+} u(\tau)$ does not exist. 
\end{example} 

\subsection{Decay of $\calJ$ and continuity at $t=0$}

Without any further knowledge on the energetic solutions, one can show that
$t\mapsto \calJ(u(t))$ in non-increasing. We emphasize that for general
energetic solutions we allow $\calJ(u^0)=\infty$, but enforce via (E) in
\eqref{eq:ES.tJ.E}  the integrability condition
$\int_0^1 \calJ(u(t))\dd t< \infty$. This will then imply the continuity $u(t)\to u(0)$
for $t\to 0^+$. See Example \ref{ex:int.J.infty} for a case with
$\int_0^1\calJ(u(t))\dd t=\infty$, where it is unclear in what sense the
limit $u^0$ is attained because the right limit $u(0^+)$ does not exist with
respect to the norm topology.  

\begin{lemma}[Decay of $\calJ$ along $u$]\label{le:RIS.J.decay}
Let $u$ be any solution of the ERIS $(X,\calE,\|\cdot\|)$. Then, we
have 
\[
\frac1{r} \|u(r)\|\geq \calJ(u(r))\geq \calJ(u(t)) \text{ \ for
} 0 < r < t .
\]
\end{lemma} 
\begin{proof} Since the dissipation is non-negative we have the energy
  estimate 
\[
e(t):=\calE(t,u(t))=t\,\calJ(u(t)) \leq r\,\calJ(u(r)) +
\int_r^t \calJ(u(\tau)) \dd \tau= e(r)+ \int_r^t \frac1\tau e(\tau)\dd \tau.
\]
Applying Gr\"onwall's estimate to $e$ we obtain $e(t)\leq (t/r) e(r)$ which
means the second estimate $\calJ(u( t))\leq \calJ(r))$.

For the first estimate we simply use stability of $u(r)$ and test
with $v=0$, namely 
\[
r\,\calJ(u(r)) \leq r\calJ(0)+\|0{-}u(r)\| = 0 + \|u(r)\|.
\]
This gives the first estimate in the assertion. 
\end{proof}

With this we now show the continuity of the energetic solutions $u$
as defined in \eqref{eq:ES.tJ}.

\begin{proposition}[Continuity at $t=0$]
\label{pr:Contin.t0}
The $u:{[0,\infty[} \to X$ be an energetic solution in the sense of
\eqref{eq:ES.tJ}, i.e.\ in particular $\int_0^1 \calJ(u(t))\dd t <\infty$. 

Then, $u$ has a right $u(0^+):=\lim_{\tau\to 0^+}
u(\tau)$ in the norm sense, and we have 
\begin{equation}
  \label{eq:VarBound}
  \mafo{Var}_{\|\cdot\|}(u;[0,t]) \leq \int_0^t\calJ(u(\tau))\dd \tau <\infty \
  \text{ for all } t>0.
\end{equation}
\end{proposition}
\begin{proof} Using $\calJ(u(t))\geq 0$, the energy balance (E) gives
\[
\mafo{Var}_{\|\cdot\|}(u;[r,t]) \leq r\calJ(u(r)) + \int_r^t \calJ(u(\tau))\dd
\tau \leq \int_0^t  \calJ(u(\tau))\dd
\tau < \infty,
\]
where we used the monotonicity $\calJ(u(\tau))\geq \calJ(u(r))$ for $\tau \in
{]0,r[}$. Since $0<r<t$ were arbitrary, we have 
\[
\mafo{Var}_{\|\cdot\|}(u;{]0,t]})= \lim_{r\to
    0^+}\mafo{Var}_{\|\cdot\|}(u;[r,t]) 
\leq   \int_0^t  \calJ(u(\tau))\dd
\tau < \infty , 
\]
which implies that $u(0^+)$ exists.  Now using the initial condition
$u(0)=u(0^+)$ we have
$\mafo{Var}_{\|\cdot\|}(u;{]0,t]}) = \mafo{Var}_{\|\cdot\|}(u;[0,t])$ and the
result is established.
\end{proof}

\begin{remark}[Attainment of the initial condition]
\label{rm:Attainment} It would be highly desirable to define energetic
solutions also for cases where $\mafo{Var}_{\|\cdot\|}(u;[0,1])=0$. One
possible way would be to replace the energy balance (E) by the corresponding
balance on subintervals $[r,t]$ as follows. 
\begin{align*}
  \wt{\text{(E)}}\quad &
  t\calJ(u(t)) + \mafo{Var}_{\|\cdot\|} (u;[r,t]) = r\calJ(u(r))+\int_r^t
  \calJ(u(\tau)) \dd \tau \quad \text{ for } 0 \lneqq r<t.
\end{align*}
However, we still need a relation to couple the solution $u:{]0,\infty[}\to X$
to its initial condition $u^0$, which could be done by defining the variational
interpolants  
\[
\wt u(t) = \mafo{arg\,min}\bigset{\| \wt u -u^0\| + t\calJ(\wt u)}{ \wt u \in X}
\]
and asking 
\begin{align*}
  \wt{\text{(I)}}\quad &
   \| u(t)-\wt u(t)\| \to 0 \quad \text{for } t\to 0^+.
\end{align*}
It is an open question whether this option provides a good definition leading
to existence and uniqueness. 
\end{remark}

\subsection{Existence theory for ERIS}
\label{su:ExistERIS}

A standard method of constructing energetic solutions is the method of
incremental minimization (also known as minimizing movement scheme). 
Choosing a time step $h>0$ we use the discrete times $t_k:=kh$ and define the
approximate solutions via 
\[
u^k_h :=\mafo{arg\,min}\bigset{\|u {-} u^{k-1}_h \| + \calE(kh, u)\|}{ u \in X},
\]
where $u^0_h=u^0$. By convexity and lower semicontinuity of $\calE(kh,\cdot)
= kh\calJ(\cdot)$ and the strict convexity of the norm, we obtain a unique
minimizer in each step and can thus construct the piecewise constant
interpolant 
\[
\ol u_h:{[0,\infty[}\to X \quad \text{with }\ol u_h(0)=u^0_k \text{ and } 
\ol u_h(t)=u^k_h \text{ for }t\in {]kh{-}h,kh]}.
\]
It is then standard to apply a Banach-space valued version of Helly's selection
principle to obtain a weakly convergent subsequence with a limit function $u$ 
and to derive an upper energy estimate, i.e.\ (E) on $[0,T]$ but with
``$\leq$'' instead of ``$=$'', see the general references \cite{MiThLe02VFRI,
  MieThe04RIHM, Miel05ERIS}.  

The major difficulty in concluding the proof is to show that the limit
function $u$ still satisfies the stability condition (S). Lemma
\ref{le:S-closed} guarantees strong closedness, while only weak convergence can
be inferred. Since the sets $\calS(t)$ of stable states are typically not convex
(see \eqref{eq:S.notcvx} or an example), they are also not weakly closed. 

To generate the missing strong convergence, the usual approach is to 
assume that the sublevels of $\calJ$ are compact in $X$ and that
$\calJ(u^0)<\infty$. Then, all approximations $u^k_h$ lie in such a compact set
and weak convergence turns into strong convergence, and existence follows by
the standard arguments as given in the above references. 

For completeness we state the following existence result, where the
compactness of the sublevels of $\calJ$ is slightly weakened by exploiting the 
a priori bound $\| u^k_h\| \leq \|u^0\|$. We emphasize that our main result
stated in Theorem \ref{th:ExiUni} does not impose any compactness
assumption. Moreover, it provides uniqueness, which cannot be derived directly
from the theory of energetic solutions. Thus, the results of the following
proposition are all contained in Theorem \ref{th:ExiUni}, but here we use the 
standard theory only, not relying on the equivalence to the GS
$(X,\calJ,\frac12\|\cdot\|^2)$. 

\begin{proposition}[Existence theory using compactness]\label{pr:Exist.Cmpt}
Consider the ERIS\\ $(X,\calE,\|\cdot\|)$ with $\calE(t,u)=t\calJ(u)$ where
$\calJ:X\to [0,\infty]$ is lower semicontinuous, convex, and 1-homogeneous. Impose additionally,
that the functional 
\[
\calG: X\to [0,\infty]; \ u \mapsto \calJ(u)+ \|u\|
\]
has compact sublevels. Then for all $u^0\in X$ with $\calJ(u^0)<\infty$ there
exists a energetic solution $u$ in the sense of \eqref{eq:ES.tJ} with
$u(0)=u^0$.  Moreover, this solution satisfies 
\[
\calJ(u(r))\geq \calJ(u(t)) \ \text{ and } \ \|u(r)\|\geq \|u(t)\| \quad 
\text{for } 0 \leq r< t. 
\]
\end{proposition}
\begin{proof} The only non-trivial part of the proof is to show that the
  approximations $u^k_h$ satisfy an a priori bound $\calG(u^k_h)\leq C$.
If this is done then, the standard existence 
theory (see e.g.\ \cite[Thm.\,6.3(2)]{MieThe04RIHM}) applies. 

To provide the bound on $\calG$ we analyze the incremental problems  in a
little more detail. Since $u^k_h$ is a minimizer, we have 
\[
kh\calJ(u^k_h)+ \|u^k_h{-}u^{k-1}_h\|  \leq kh\calJ(u^{k-1}_h)+ 0,
\] 
which implies $\calJ(u^k_h) \leq \calJ(u^{k-1}_h)$. 

We also claim $\|u^k_h\|\leq \|u^{k-1}_h\|$. To see this, we may restrict to
the case $u^k_h \neq u^{k-1}_h$, since otherwise the inequality holds
trivially. Then, the Euler-Lagrange equation reads 
\[
0 \in \frac1{\|u^k_h {-} u^{k-1}_h\|}\,E\big( u^k_h {-} u^{k-1}_h\big) + 
kh\, \pl\calJ (u^k_h).
\]
Testing this equation with $u^k_h$ and using that $\langle \eta ,u\rangle =
\calJ(u)\geq 0$ for all $\eta \in \pl\calJ(u)$ we conclude
\[
\|u^k_h\|^2 - \big(u^{k-1}_h\big|u^k_h\big) = \langle E\big( u^k_h {-}
u^{k-1}_h\big), u^k_h\rangle = -\|u^k_h {-} u^{k-1}_h\| \,kh\, \calJ(u^k_h)  \leq 0. 
\]
This implies $\| u^k_h\| \leq \| u^{k-1}_h\|$ as desired. 

Together, we obtain the monotonicity $\calG(u^k_h)\leq \calG(u^0)$ and conclude
that all values  $\ol u_h(t)$ of the approximating functions lie in the compact
set $K:=\bigset{u\in X}{ \calG(u)\leq \calG(u^0)}$. Thus existence follows. 

The monotonicity of $\calJ$ follows from Lemma \ref{le:RIS.J.decay}. 
Moreover, for the approximation functions $\ol u_h$ we have $\| \ol
u_h(r)\|\geq \| \ol u_h(t)\|$ for $0 \leq r< t$. Because the weak convergence
in the compact set $K$ is turned into strong convergence, this inequality
survives for the limit function as well. 
\end{proof}

It is an open question how to show the monotonicity of $t\mapsto \|u(t)\|$
directly for all energetic solutions.

\subsection{Advanced properties of energetic solutions}
\label{su:AdvPropES}
The next result characterizes jumps and shows that along a jump the solution
can be modified with any value on the straight line connecting the left and the
right limit. The result is essentially contained in \cite{MieThe04RIHM}, but we
provide a full proof for the present special case.

To state the result we introduce the notation of left and right
limits $u(t^\pm)$ of $u:{[0,\infty[} \to  X$, which exist since
$\DISS(u,[0,T])$ is finite for all $T>0$. We set 
\[
u(t^-)= \lim_{\tau\to t^-} u(\tau) \quad \text{ and } \quad
u(t^+)= \lim_{\tau\to t^+} u(\tau),
\]
and use a corresponding notation for the dissipation, namely 
\[
\DISS(u,{[t_1,t_2[})=\lim_{\tau\to t^-_2}
\DISS(u,[t_1,\tau]) = \DISS(u,[t_1,t_2]) - \|
u(t_2){-}u(t_2^-)\|.
\] 
In the following result the assertions (i) to (iii) hold even without
convexity, see \cite[Lem.\,2.1.13]{MieRou15RIST}. For (iv) convexity is needed,
but not the 1-homogeneity. 

\begin{proposition}[Jumps in ERIS]\label{pr:RIS.jumps}
Consider a solution $u:{[0,\infty[} \to X$ of ERIS and a time $t>0$ such
$u$ is not continuous at $t$ (i.e.\ not all of the three values
$u(t^-)$, $u(t)$, and $u(t^+)$ are the same). Then, we have the
relations 
\begin{equation}
  \label{eq:Jumps}
  \begin{aligned} 
\text{(i) }  &\calE(t,u(t^+))+\|u(t^+){-}u(t^-)\|= \calE(t,u(t^-)), \\
\text{(ii) }  &\calE(t,u(t))+\|u(t){-}u(t^-)\|= \calE(t,u(t^-)), \\
\text{(iii) }  &\exists\, \theta_* \in [0,1]: \quad u(t)=(1{-}\theta_*)
u(t^-) + \theta_* u(t^+),\\
\text{(iv) }&\forall\,\theta\in[0,1]:\  \calE\big(t,(1{-}\theta)
u(t^-) {+} \theta u(t^+)\big)= (1{-}\theta)\calE(t,u(t^-)) + \theta
\calE(t,u(t^+))  .
  \end{aligned} 
\end{equation}
Indeed, if we modify $u$ at the time $t$ by replacing $\theta_*$ in (iii) by
any other $\theta\in [0,1]$, we still have a solution of ERIS, in particular
$(1{-}\theta) u(t^-) {+} \theta u(t^+) \in \calS(t)$.
\end{proposition}
\begin{proof} The upper energy estimates 
\[
\calE(t,u(t^+))+\|u(t^+){-}u(t)\|\leq \calE(t,u(t)) \text{ and } 
\calE(t,u(t))+\|u(t){-}u(t^-)\|\leq \calE(t,u(t^-)) 
\]
follow from the energy balance on $[t,t_2]$ and $[t_1,t]$ and taking
the limits $t_2\to t^+$ and $t_1\to t^-$, respectively. The
lower estimates follow from the stability of $u(t)$ and $u(t^-)$,
respectively. For the latter stability use Lemma \ref{le:S-closed} and
$\calS(t_n) \ni u(t_n)\to u(t^-)$ for $t_n\to t^-$. 

Having the two identities we obtain by summing 
\begin{align*}
\calE(t,u(t^-)) &= \calE(t,u(t^+)) + \|u(t^-){-}u(t)\|+
\|u(t){-}u(t^+)\|\\ &\geq \calE(t,u(t^+))+\|u(t^+){-}u(t^-)\| \geq
\calE(t,u(t^-)),  
\end{align*}
where the last estimate follows from the stability of $u(t^-)$. We conclude
$\|u(t^-){-}u(t)\|+ \|u(t){-}u(t^+)\| =\|u(t^+){-}u(t^-)\|$, which implies
(iii), since we are in a Hilbert space.

To establish (iv) we  use the abbreviation $u_\theta:=(1{-}\theta)
u(t^-) + \theta u(t^+)$. On the one hand, the convexity of $\calJ$ and the
established identity give the upper estimate
\[
\calE(t,u_\theta) \leq \calE(t,u_0) - \theta \| u_1{-}u_0\|.
\]
On the other hand, the stability of $u_0=u(t^-)$ gives the lower estimate
\[
\calE(t,u_\theta) \geq  \calE(t,u_0)  - \| u_\theta{-}u_0\|
=\calE(t,u_0) - \theta \| u_1{-}u_0\|.
\]
The last two estimates imply (iv). The stability of $u_\theta$ now
follows from $u_0 \in \calS(t)$, namely 
\[
\calE(t,u_\theta)=\calE(t,u_0)- \| u_\theta{-}u_0\| \leq \calE(t,\wt
u) + \|\wt u{-}u_0\|-\|u_\theta{-}u_0\| \leq \calE(t,u_0) + \|\wt u
{-}u_\theta\|. 
\] 
This proves the result. 
\end{proof}

\subsection{Time-dependent dissipation}
\label{su:TimeDepDiss}

Instead of looking at the time-dependent energy $\calE(t,u)=t\calJ(u)$
and the time-independent dissipation $\|\dot\|$, we may multiply the
equation by $1/t$ to obtain the ERIS $(X,\calJ(u),\frac1t\|\cdot\|)$,
where now the time-dependence is in the dissipation functional
\[
\calR(t,\dot u)=\frac1t\|\dot u\|. 
\]
The stability sets $\calS(t)$ are still the same as well as the
differential form of the (formal) power balance: 
\[
\frac{\rmd}{\rmd t} \calE(t,u(t)) + \|\dot u(t)\| = \pl_t\calE(t,u(t))
\quad \Longleftrightarrow \quad 
\frac{\rmd}{\rmd t} \calJ(u(t)) + \frac1t\|\dot u(t)\| =0. 
\]
Thus, the rigorously formulated energy balance reads 
\begin{equation}
  \label{eq:EB1/t}
  \calJ(u(t)) + \int_{[s,t]} \frac1\tau \|\rmd u(\tau)\| = \calJ(u(s))
\text{ for } 0 < s < t.
\end{equation}
For a general continuous function $\phi:{]0,\infty[}\to \R$ and $0<s<t$, we can define the
weighted variation via
\begin{equation}
  \label{eq:WeightVar}
 \begin{aligned} \int_{[s,t]} \phi(\tau)\ \|\rmd u(\tau)\|:= \sup\Big\{& 
\sum_{j=1}^N \min_{\tau\in [t_{j-1},t_j]} \phi(\tau)
 \ \|u(t_j){-}u(t_{j-1})\| \\
& \qquad \Big|\;  N\in \N, \ s=t_0{-}t_1{<}\cdots{<}t_N=t\Big\},
\end{aligned}
\end{equation}
see also \cite[App.\,B.5]{MieRou15RIST}.

\section{From GS to ERIS} 
\label{se:GS2ERIS}

We now show that the solutions $w(s)=\mfT_s(w(0))$ give rise to
energetic solutions for the ERIS $(X,\calE,\|\cdot\|)$. 
For this we define the mapping $S:{[0,\infty[} \to {[0,\infty[} $ via the
relation 
\[
S(0)= 0 \quad \text{ and } S(t) = \min\bigset{ s \geq 0 }{  \|w'_+(s)\|
  \leq 1/t  } \text{ for }t>0. 
\]
Note that $g_w:s\mapsto  \|w'_+(s)\|$ is non-increasing and continuous from the right,
hence it is lower semicontinuous and the minimum is really attained.
In particular, we have $g_w(S(t))= \|w'_+(S(t))\| \leq 1/t$ by construction.  
As a consequence, $S$ is non-decreasing and continuous
from the left. Relation \eqref{eq:Veloc} provides the upper bound  $S(t)
\leq t\,\sqrt2\,\|w(0)\|$.

If $g_w$ has a jump at $s_*>0$ with
\[ 
a_*= g_w(s_*^-):=\lim_{s\to s_*^-} g_w(s) \gneqq g_w(s_*)=b_*,
\]
then $S$ has a plateau with $S(t)=s_*$ for $t\in {]1/a_*,1/b_*]}$. 
Vice versa, if $g_w$ has a plateau ${[s_1,s_2[}$ with $g_w(s)=a_*>0$,
then $S$ has a jump in $t_*:=1/a_*$ with $s_1=S(t_*)=\lim_{t\to
  t_*^-} S(t)$  and $s_2=\lim_{t\to t^+_*} S(t)$.  If $g_w$ has a
plateau with value $a_*=0$, then the plateau is ${[s_1,\infty[}$, and
$S$ remains bounded by $s_1$. 

We now define the function 
\begin{equation}
  \label{eq:def.u}
  u:\left\{\ba{ccc}{[0,\infty[} &\to &X,\\ t &\mapsto &w(S(t)).  \ea \right. 
\end{equation}
By construction, the function $u$ is continuous from the left, because
$w:{[0,\infty[} \to X$ is continuous and $S:{[0,\infty[} \to {[0,\infty[}$ is
continuous from the left.

The next result is crucial for connecting the gradient system
$(X,\calJ,\frac12\|\cdot\|^2)$ with the RIS $(X,\calE,\|\cdot\|)$. 

\begin{proposition}[Stability]\label{pr:w.stable}
Consider a solution $w$ of $(X,\calJ,\frac12\|\cdot\|^2)$, then
\begin{equation}
  \label{eq:w-stable}
  \forall\ s>0:\quad w(s) \in \calS\big(1/\|w'_+(s)\|\big).
\end{equation}
\end{proposition}
\begin{proof}
For $s>0$ the right derivative $w'_+(s)\in X$ exists and $0\in Ew'_+(s) +
\pl\calJ(w(s))$. Then, the convexity of $\calJ$ gives the estimate
\begin{align*}
\calJ(v)\geq \calJ(w(s))+ \langle Ew'_+(s),v{-}w(s)\rangle \geq
\calJ(w(s)) - \|w'_+(s)\| \|v{-}w(s)\|.
\end{align*}
Together with $\calE(t,u)=t \calJ(u)$ this means $w(s) \in
\calS(t)$ for $t=1/\|w'_+(s)\|$. 
\end{proof} 

We are now ready to establish our main existence result for the ERIS,
which is obtained as a consequence of the existence result for the
gradient system and the corresponding reparametrization. We emphasize that we
do not assume any type of compactness.

\begin{theorem}[From GS to ERIS]\label{th:GF2RIS}
Let $w:{[0,\infty[} \to X$ be a solution of the gradient system
$(X,\calJ,\frac12\|\cdot\|^2)$ 
with $\calJ(w(0))<\infty$,
then the function $u:{[0,\infty[} \to X$
defined in  \eqref{eq:def.u} is an energetic solution for the ERIS
$(X,\calE,\| \cdot\|)$ in the sense of \eqref{eq:ES.tJ}. 
\end{theorem}
\begin{proof} For $w(0)=0$ we have $w\equiv 0$ and hence $u\equiv 0$, which is
  trivially an energetic solution of ERIS. Thus, we now assume $w(0)\neq 0$.

\emph{Stability (S):\/} For $t=0$ the stability $u(0)=w(0)\in
\calS(0)=X$ is trivial. For $t>0$ we have $s=S(t)>0$ and conclude 
\[
u(t)=w(S(t)) \in S\big(1/g_w(S(t))\big) \subset \calS(t),
\]
where we used $g_w(S(t))\leq 1/t$. Hence, \eqref{eq:ES.tJ.S} is established.

\emph{Energy balance (E):\/} According to the general theory of RIS, it
is sufficient to establish an upper energy estimate, since the lower
estimate is a consequence of the stability, cf.\
\cite[Prop.\,2.1.23]{MieRou15RIST} or \cite{MieThe04RIHM}. 

In principle the energy balance follows from the energy-dissipation balance
\eqref{eq:EDB-GF} for the gradient system, where we would like to use the time
reparametrization $s=S(t)$ or $t=g_w(s)=1/\|w'_+(s)\|$ giving formally
$\|w'(s)\|\dd s = \dd t$. However, because of jumps we have to be more careful
and estimate the dissipation explicitly.  We have
\[
\DISS(u,[r,t]) = \sup\Bigset{\sum_{j=1}^N
  \|u(t_j){-}u(t_{j-1})\|}{ N\in \N, \ r= t_0<t_1 < \cdots <
  t_N= t}.
\]
We choose a finite partition $(t_j)_{j=0,1,..,N}$ of $[r,t]$ such that
$\DISS(u,[t,T])$ is approximated  up to an error smaller than $\eps$ 
and set $s_j=S(t_j)$. 
The monotonicity of $g_w$ and \eqref{eq:EDB-GF} yield
\begin{align*}
\calJ(w(s_{j-1})) &= \calJ(w(s_j))+\int_{s_{j-1}}^{s_j} \| w_+'(s)\|^2 \dd s
\geq \calJ(w(s_j))+\|w'(s_j)\| \int_{s_{j-1}}^{s_j} \| w_+'(s)\| \dd s\\
&\geq \calJ(w(s_j))+\frac1{t_j} \| w(s_j){-}w(s_{j-1})\|.
\end{align*}
In terms of the ERIS $(X,\calE,\|\cdot\|)$ and the function $u$, this
means 
\[
\calE(t_j,u(t_j)) + \|u(t_j){-}u(t_{j-1})\| \leq
\calE(t_{j-1},u(t_{j-1})) + \int_{t_{j-1}}^{t_j} \pl_\tau
\calE(\tau,u(t_{j-1})) \dd \tau.
\]
Summing over $j\in \{1,..,N\}$  we obtain 
\[
  \calE(t,u(t))+ \mafo{Var}_{\|\cdot\|}(u;[r,t])\leq \eps + \calE(r,u(r)) +
  \sum_{j=1}^N \int_{t_{j-1}}^{t_j} \pl_\tau \calE(\tau,u(t_{j-1})) \dd \tau.
\]   
Taking $\eps$ and the fineness of the partition to $0$ simultaneously,
we obtain the desired energy estimate \eqref{eq:ES.tJ.E} 
on all intervals $[r,t]$ with $0<r<t$. 

\emph{Initial condition (I):\/} Using $\calJ(w(0))=\calJ(u(0)) <\infty$ and the
monotonicity of $\calJ$, we obtain $\int_0^1 \calJ(u(t))\dd t \leq 
\calJ(u(0))<\infty$.  Thus, Proposition
\ref{pr:Contin.t0} provides the desired continuity, and \eqref{eq:ES.tJ} is
established as well.
\end{proof}

\section{From ERIS to GS}
\label{se:RIS2GF}

Here we show that every energetic solution $u:{[0,\infty[} \to X$ for the
ERIS $(X,\calE,\|\cdot\|)$ gives rise to a solution $w:{[0,\infty[} \to X$
for the GS $(X,\calJ,\frac12\|\cdot\|^2)$. We do this by
reparametrization. However, at jumps we need to fill in pieces, which
can be done in a piecewise affine manner. 

Affine interpolations are defined for each function
$g:{[0,\infty[} \to V$ having left and right limits $g(t^\pm)$ for all $t$,
where $g(0^-):=g(0)$.  The interpolant reads
\[
g_\theta(t)= (1{-}\theta) g(t^-) + \theta g(t^+), \text{ where
}\theta\in [0,1] \text{ and } t\geq 0, 
\]

The time reparametrization is given in terms of the left-continuous
function
\begin{align*}
\wh s(t)&=\int_0^t \!\tau\|\rmd u(\tau)\|
\\
&:= \sup\Bigset{\sum_{j=1}^N t_{j-1}
  \| u(t_j){-}u(t_{j-1})\| }{N\in \N, \ 0\leq t_0{<}t_1{<}\cdots {<}t_N\leq t}.  
\end{align*}
An inverse of this function is given by 
\[
 \wh t(s) :=
\inf\set{ t\geq 0 }{ \wh s (t)\geq s } . 
\] 
Hence, $\wh s$ will have a jump at $t_*$, if $u$ has a jump at $t_*$, more
precisely $\wh s(t^+_*)-\wh s(t_*)= t_*\|u(t^+_*){-}u(t_*)\|$. 
In contrast, $\wh t$ will have a plateau, namely $\wh t(s)= t_*$ for
$s\in {[\wh s(t_*),\wh s(t^+_*)[}$.
Moreover, if $\wh s$ has a plateau for ${]t_1,t_2]}$ with value $s_*$, 
then $\wh t$ has a jump at $s_*$.  

For a given energetic solution $u:{[0,\infty[} \to X$, we define the
function 
\begin{equation}
  \label{eq:w.from.u}
  w(s) =  u_\theta(\wh t(s)) \text{ for } s= \wh \sigma_\theta(s), \quad 
\text{where } \wh\sigma(s):=\wh s(\wh t(s)).
\end{equation}
Note that $\wh \sigma(s)=\min\bigset{\wh s(t)}{t\geq 0,\ \wh s(t)\leq s }
\leq s$ and $\wh\sigma(s)< s$ only in regions ${]s_1,s_2[}$ which are
not covered by the range of $\wh s $, i.e.\ there exists $t_*$ such
that $\wh s(t_*)\leq s_1 < s_2 \leq \wh s(t_*^+)$ and $\wh t(s)=t_*$.   

\begin{theorem}[From energetic solutions to gradient-flow
  solutions] \label{th:RIS2GF} 
If $u:{[0,\infty[} \to X$ is an energetic solution for
$(X,\calE,\|\cdot\|)$, then the function $w:{[0,\infty[} \to X$ defined in 
\eqref{eq:w.from.u} is a solution for the gradient system
$(X,\calJ,\frac12\|\cdot\|^2)$, i.e.\ it satisfies \eqref{eq:GFE}.  
\end{theorem}
\begin{proof} \emph{Step 1:}\/  We first show that $w$ lies in
  $\rmH^1_\text{loc}({]0,\infty[};X)$. 
For this, we need to estimate $\frac1{s_2{-}s_1} \| w(s_2){-}w(s_1)\|^2$ as
follows. We have  $s_j =(1{-}\theta_j) \wh s(t_j)+ \theta_j \wh
s(t_j^+)$ for suitable $\theta_j \in [0,1]$. With this choice and the
definition of $\wh s$ based on the variation of $u$ we obtain
\begin{equation}
  \label{eq:s2-s1}
  \begin{aligned}
s_2{-}s_1&= \theta_2\big(\wh s(t_2^+)- \wh s(t_2) \big) + \big( \wh s(t_2) -
\wh s(t_1^+) \big) + (1{-}\theta_1) \big(\wh s(t_1^+) - \wh s(t_1)
\big)\\
& \geq \theta_2t_2 \| u(t_2^+){-}u(t_2)\| + t_1\DISS(u,{]t_1,t_2[}) +
  (1{-}\theta_1)t_1\| u(t_1^+){-}u(t_1)\| \\
& \geq t_1\| \, u_{\theta_2}(t_2) - u_{\theta_1}(t_1)\,\| \ = \ t_1\|
w(s_2){-}w(s_1)\|.  
\end{aligned}
\end{equation}
We conclude that for all $0< s_1<s_2$ we have 
\[
 \frac1{s_2{-}s_1} {\| w(s_2){-}w(s_1)\|^2}  
  = \frac{\| w(s_2){-}w(s_1)\| }{s_2{-}s_1} \: \|
                   u_{\theta_2}(t_2) {-} u_{\theta_1}(t_1) \| 
  \leq \frac1{t_1}\,\| u_{\theta_2}(t_2) {-} u_{\theta_1}(t_1) \|  .
\]
Hence, for any partition $0<r=s_0< s_1 < \cdots <s_{N-1}<s_N=s$ of
$[r,s]$ we
obtain 
\[
\sum_{j=1}^N \frac{ \|w(s_j){-}w(s_{j-1}\|^2}{s_j{-}s_{j-1}} \leq 
\sum_{j=1}^N \frac1{t_{j-1}} \|
u_{\theta_j}(t_j){-}u_{\theta_{j-1}}(t_{j-1})\| 
\leq \frac1{t_0}\DISS\big(u, [\wt t(r), \wt t(s^+)] \big)< \infty. 
\]
Since the partition was arbitrary and since every energetic solution
has bounded variation on all intervals compactly contained in
${]0,\infty[}$, we conclude $w|_{[r,s]} \in \rmH^1([r,s];X)$,
which is the desired result.\medskip 

\emph{Step 2:}\/ More precisely, for $s_j$ with $\wh
\sigma(s_j)=s_j$  we have 
\[
\int_{s_1}^{s_2} \|w'(s)\|^2 \dd s \leq  \int_{\wt t(s_1)}^{\wt t(s_2)}
\frac 1t \|\rmd u(t)\|= \calJ\big(u(\wt t(s_1))\big) {-} 
\calJ\big(u(\wt t(s_2))\big) = \calJ(w(s_1)){-} \calJ(w(s_2)) ,
\]
where we used the energy balance \eqref{eq:EB1/t} for the time-dependent
dissipation model.

Moreover, dividing \eqref{eq:s2-s1} by $(s_2{-}s_1)t_1$ we may pass to the limit
$s_2\to s^+_1$ for  almost all $s_1=s>0$  and obtain 
\[
\frac1{t_1}= \frac1{\wh t(s)} = \| w'_+(s)\|. 
\]   

\emph{Step 3:}\/ We now want to show that $w$ solves the
gradient-flow equation \eqref{eq:GFE}, i.e.\ $0\in Ew'(s) +
\pl\calJ(w(s))$ for a.a.\ $s>0$.  For this we use the stability of
$u$. Indeed by the construction of the interpolant $u_\theta(t)$ we
have the stability of $u_\theta(t)\in \calS(t)$, see Proposition
\ref{pr:RIS.jumps}. Since $u_\theta(t)$ is a minimizer of $\wt u
\mapsto t\calJ(\wt u) + \|\wt u{-}u_\theta(t)\|$, we know that 
$\pl\calJ(u_\theta(t))$ is
nonempty and that $\| \pl^0\calJ(u_\theta(t))\|_* \leq 1/t$. 
Translating this to the variable $s=\wh s (t)$, we obtain 
\begin{equation}
  \label{eq:eta.w'}
  \forall\,s>0 : \quad
\|\eta(s)\|\leq 1/\wh t(s)= \|w'(s)\| \quad \text{with } \eta(s)=\pl^0\calJ(w(s)) , 
\end{equation} 
where the last relation follows with Step 2. 

Since $\calJ$ is convex and $w|_{[s_1,s_2]}\in \rmH^1([s_1,s_2];X)$
with $[s_1,s_2]\ni s\mapsto \calJ(w(s))$ bounded and decaying we can
apply the chain rule (see \cite[Lem.\,3.3]{Brez73OMMS}) and obtain 
\[
\frac{\rmd}{\rmd s} \calJ(w(s)) = \langle \eta(s),w'(s)\rangle \quad
\text{ a.e.} 
\]
Integration gives the first identity in the following relations, 
and exploiting \eqref{eq:eta.w'} and the energy estimate
from Step 2 yields 
\begin{align*}
\calJ(w(s_2))-\calJ(w(s_1)) &= 
\int_{s_1}^{s_2}\langle \eta(s),w'(s)\rangle \dd s 
\geq \int_{s_1}^{s_2}{-} \| \eta(s)\|\,\| w'(s)\| \dd s \\
&  \overset{\text{\eqref{eq:eta.w'}}}\geq 
- \int_{s_1}^{s_2}\| w'(s)\|^2 \dd s  
 \overset{\text{Step 2}}\geq   \calJ(w(s_2))-\calJ(w(s_1)) .
\end{align*}  
Thus, we conclude that all ``$\geq$'' must be equalities. The first estimate
shows $\eta(s)= -\alpha(s) Ew'(s)$ for some $\alpha(s)\geq 0$. By
\eqref{eq:eta.w'} we know $\alpha(s) \leq 1$, but then the second
estimate gives $\alpha(s)=1$ for a.e.\ $s>0$. Thus, $\eta(s)=-Ew'(s)$
for a.e.\ $s>0$ gives the desired equation $0\in Ew'(s) +
\pl\calJ(w(s))$. 
\end{proof}

With the available link from the ERIS to the GS, we are now able to conclude
the proof of our main theorem on the ERIS by exploiting the existence and
uniqueness results for the GS in Section \ref{se:GS}, which do not need any
compactness assumption.
\\[0.4em]
\begin{proof}[Proof of Theorem \ref{th:ExiUni}] \hspace*{-1em}
Existence of solution for the ERIS follows from 
Theorem \ref{th:GF2RIS}, which shows that suitably reparametrizing
the solutions of the GS leads to energetic solutionsin the sense
of \eqref{eq:ES.tJ} for the ERIS. 

The uniqueness of energetic solutions follows using Theorem
\ref{th:RIS2GF}, since every energetic solution generates a solution
of the gradient system. Since the latter is unique, the uniqueness of
energetic solutions follows if we choose the unique left-continuous variant
and neglect the freedom to choose the
value $u_\theta(t)=(1{-}\theta)u(t^-)+\theta u(t^+)$ of an energetic solution   
at an jump time $t$. 
\end{proof}

In our main theorem, we have restricted the existence and uniqueness result for
the ERIS to initial values $u^0$ with $\calJ(u^0)<\infty$. 
For the GS the existence and uniqueness result extends to all 
initial conditions $u^0\in X$ because dom$(\calJ)$ is assumed to be
dense. However, there is a subtle issue about the rate-independent rescaling
$u(t)=w(S(t))$ which may lead to $\mafo{Var}_{\|\cdot\|}(u;[0,1])=\infty$ when
doing the corresponding reparametrization. It remains
an open problem to provide an intrinsic formulation of energetic solutions and
their attainment of the initial condition in the case of infinite variation
near $t=0$, see Remark \ref{rm:Attainment}.

\paragraph*{Acknowledgments.}
The research was partially supported by Deutsche Forschungsgemeinschaft (DFG)
via the Collaborative Research Center SFB\,910 ``Control of self-organizing
nonlinear systems'' (project number 163436311), subproject A5 ``Pattern
formation in coupled parabolic systems''. The author is grateful to Martin
Burger for stimulating discussions.

\footnotesize

\begin{thebibliography}{11}\itemsep0.1em

\bibitem[ACM04]{AnCaMa04PQEM}
{\scshape F.~{Andreu-Vaillo}, V.~Caselles, {\upshape and} J.~M.~Maz\'on}.
\newblock {\em Parabolic Quasilinear Equations Minimizing Linear Growth
  Functionals}.
\newblock Birkh\"auser-Verlag, 2004.
\newblock xiv+340 pp.

\bibitem[AGS05]{AmGiSa05GFMS}
{\scshape L.~Ambrosio, N.~Gigli, {\upshape and} G.~Savar{\'e}}.
\newblock {\em Gradient flows in metric spaces and in the space of probability
  measures}.
\newblock Lectures in Mathematics ETH Z\"urich. Birkh\"auser Verlag, Basel,
  2005.

\bibitem[BCN02]{BeCaNo02TVFR}
{\scshape G.~Belletini, V.~Caselles, {\upshape and} M.~Novaga}.
\newblock The total variation flow in {$\mathbb R^N$}.
\newblock {\em J. Diff. Eqns.}, 184, 475--525, 2002.

\bibitem[BDM15]{BoDuMa15TDVA}
{\scshape V.~B{\"o}gelein, F.~Duzaar, {\upshape and} P.~Marcellini}.
\newblock A time dependent variational approach to image restoration.
\newblock {\em SIAM J. Imag. Sci.}, 8(2), 968--1006, 2015.

\bibitem[BKS04]{BrKrSc04UEQI}
{\scshape M.~Brokate, P.~Krej{\v{c}}{\'{\i}}, {\upshape and} H.~Schnabel}.
\newblock On uniqueness in evolution quasivariational inequalities.
\newblock {\em J. Convex Anal.}, 11, 111--130, 2004.

\bibitem[BoF12]{BonFig12TVFS}
{\scshape M.~Bonforte {\upshape and} A.~Figalli}.
\newblock Total variation flow and sign fast diffusion in one dimension.
\newblock {\em J. Diff. Eqns.}, 252(8), 4455--4480, 2012.

\bibitem[Br{\'e}73]{Brez73OMMS}
{\scshape H.~Br{\'e}zis}.
\newblock {\em Op{\'e}rateurs maximaux monotones et semi-groupes de
  contractions dans les espaces de {H}ilbert}.
\newblock North-Holland Publishing Co., Amsterdam, 1973.

\bibitem[HNV19]{HeNeVa19?SHLC}
{\scshape M.~Heida, S.~Neukamm, {\upshape and} M.~Varga}.
\newblock Stochastic homogenization of {$\Lambda$}-convex gradient flows.
\newblock {\em arXiv:1905.02562}, 2019.

\bibitem[KMR13]{KiMuRy13ACST}
{\scshape K.~Kielak, P.~B.~Mucha, {\upshape and} P.~Rybka}.
\newblock Almost classical solutions to the total variation flow.
\newblock {\em J. Evol. Eqn.}, 13, 21--49, 2013.

\bibitem[Mie05]{Miel05ERIS}
{\scshape A.~Mielke}.
\newblock Evolution in rate-independent systems ({C}h.~6).
\newblock In C.~Dafermos {\upshape and} E.~Feireisl, editors, {\em Handbook of
  Differential Equations, Evolutionary Equations, vol.~2}, pages 461--559.
  Elsevier B.V., Amsterdam, 2005.

\bibitem[Mie08]{Miel08CIME}
{\scshape A.~Mielke}.
\newblock Differential, energetic and metric formulations for rate-independent
  processes.
\newblock Slides of Lecture Series given at ``C.I.M.E. Summer School on
  Nonlinear PDEs and Applications'', Cetraro, June 2008.

\bibitem[Mie16]{Miel16EGCG}
{\scshape A.~Mielke}.
\newblock On evolutionary {$\Gamma$}-convergence for gradient systems
  {(Ch.\,3)}.
\newblock In A.~Muntean, J.~Rademacher, {\upshape and} A.~Zagaris, editors,
  {\em Macroscopic and Large Scale Phenomena: Coarse Graining, Mean Field
  Limits and Ergodicity}, Lecture Notes in Applied Math. Mechanics Vol.\,3,
  pages 187--249. Springer, 2016.
\newblock Proc. of Summer School in Twente University, June 2012.

\bibitem[MiR07]{MieRos07EURC}
{\scshape A.~Mielke {\upshape and} R.~Rossi}.
\newblock Existence and uniqueness results for a class of rate-independent
  hysteresis problems.
\newblock {\em Math. Models Meth. Appl. Sci. (M$^3$AS)}, 17(1), 81--123, 2007.

\bibitem[MiR15]{MieRou15RIST}
{\scshape A.~Mielke {\upshape and} T.~Roub{\'\i}{\v{c}}ek}.
\newblock {\em Rate-Independent Systems: Theory and Application}.
\newblock Applied Mathematical Sciences, Vol.\,193. Springer New York, 2015.

\bibitem[MiT04]{MieThe04RIHM}
{\scshape A.~Mielke {\upshape and} F.~Theil}.
\newblock On rate--independent hysteresis models.
\newblock {\em Nonl. Diff. Eqns. Appl. (NoDEA)}, 11, 151--189, 2004.
\newblock (Accepted July 2001).

\bibitem[MTL02]{MiThLe02VFRI}
{\scshape A.~Mielke, F.~Theil, {\upshape and} V.~I.~Levitas}.
\newblock A variational formulation of rate--independent phase transformations
  using an extremum principle.
\newblock {\em Arch. Rational Mech. Anal.}, 162, 137--177, 2002.

\bibitem[QuM08]{QuaMar08ADCR}
{\scshape J.~Quah {\upshape and} D.~Margetis}.
\newblock Anisotropic diffusion in continuum relaxation of stepped crystal
  surfaces.
\newblock {\em J. Phys. A: Math. Theor.}, 235004, 18pp, 2008.

\bibitem[ROF92]{RuOsFa92NTVB}
{\scshape L.~I.~Rudin, S.~Osher, {\upshape and} E.~Fatemi}.
\newblock Nonlinear total variation based noise removal algorithms.
\newblock {\em Physica D}, 60, 259--268, 1992.

\end{thebibliography}

\def\cprime{$'$}

\end{document}